\newif\ifMA
\MAfalse

\ifMA\else	
  \newif\iftwoside
  \twosidefalse
\fi

\newcounter{version}
\ifMA
    \documentclass[reqno,12pt,a4paper,oneside]{amsart}
\else
  \ifnum\theversion>0
    \documentclass[reqno,12pt,a4paper,oneside,draft]{amsart}
    
  \else
    \iftwoside
      \documentclass[reqno,12pt,a4paper,twoside]{amsart}
    \else
      \documentclass[reqno,12pt,a4paper,oneside]{amsart}
    \fi
  \fi
\fi
\usepackage{amscd}

\newcommand{\Wedge}{\mathop{\text{\larger$\wedge$}}\nolimits}
\newcommand{\Otimes}{\mathop{\text{\larger$\otimes$}}\nolimits}
\newcommand{\Oplus}{\mathop{\text{\larger[2]$\oplus$}}}
\newcommand{\Cup}{\mathop{\text{\larger$\cup$}}}

\renewcommand{\Bbb}{\mathbb}
\newcommand{\abs}[1]{\lvert#1\rvert}

\newcommand{\Cn}{{\Bbb C}^n}
\newcommand{\D}{\Delta}

\newcommand{\db}{\bar\partial}
\newcommand{\Proj}{{\bf P}}
\newcommand{\Hypb}{{\bf H}}
\newcommand{\R}{{\Bbb R}}
\newcommand{\im}{\sqrt{-1}}
\newcommand{\opn}[1]{\qopname\relax o{#1}}
\newcommand{\Un}{\opn U}
\newcommand{\Clif}{\opn{CLIF}}
\newcommand{\Spin}{\opn{SPIN}}
\newcommand{\Or}{\opn{O}}
\newcommand{\SO}{\opn{SO}}
\newcommand{\GL}{\opn{GL}}
\newcommand{\SU}{\opn{SU}}
\newcommand{\E}{\opn{E}}
\newcommand{\gl}{\opn{\mathfrak{gl}\,}}

\newcommand{\Pin}{\opn{PIN}}
\newcommand{\spin}{\opn{\mathfrak{spin}\,}}
\newcommand{\un}{\opn{{\mathfrak u}\,}}

\newcommand{\C}{{\Bbb C}}

\newcommand{\PT}{{\mathcal P}}
\newcommand{\AAA}{{\mathcal Q}}
\newcommand{\FF}{{\mathcal F}}
\newcommand{\partdel}[2]{\frac{\partial#1}{\partial#2}}
\newcommand{\ol}[1]{\bar{#1}}
\newcommand{\Db}{D_{\beta}}
\newcommand{\Da}{D_{\alpha}}
\newcommand{\Dr}{D_{\gamma}}
\newcommand{\Ric}[2]{{\mathcal R}_{#1\ol{#2}}}
\newcommand{\Curv}[4]{R_{#1\ol{#2}#3\ol{#4}}}
\newcommand{\dda}{d_{\alpha}}
\newcommand{\ddb}{d_{\beta}}
\newcommand{\ddv}{d_{\gamma}}
\newcommand{\Eb}{E_{\beta}}
\newcommand{\Ea}{E_{\alpha}}
\newcommand{\Er}{E_{\gamma}}
\newcommand{\Lb}[1]{\tilde{L_{\ol{#1}}}}
\newcommand{\La}[1]{\tilde{L_{#1}}}
\newcommand{\jjj}{j(\phi)}
\newcommand{\iii}{{\ol{I}_1,\dots,\ol{I}_h}}
\newcommand{\dwa}{\widehat {dw_{ij}}}
\newcommand{\dwb}{\widehat {d{\ol w}_{ij}}}
\newcommand{\cO}{{\mathcal O}}
\newcommand{\hotimes}{{\hat\otimes}}
\newcommand{\hS}{\hat{S}}
\newcommand{\El}{\mathcal{L}}
\newcommand{\gj}{g(\Da)\jjj}
\newtheorem{thm}{Theorem}[section]

\newtheorem{lem}[thm]{Lemma}
\theoremstyle{definition}
\newtheorem{defn}[thm]{Definition}
\theoremstyle{remark}
\newtheorem*{rem}{Remark}
\theoremstyle{remark}
\newtheorem*{expl}{Example}
\numberwithin{equation}{section}

\newcommand{\Romannumeral}[1]
  {\uppercase\expandafter{\romannumeral #1}}
\newcommand{\reftag}[1]{\tag{\ref{#1}}}

\begin{document}

\title[The Penrose transform on conformally Bochner-K\"ahler manifolds]{}
\author[Yoshinari Inoue]{}

{\quad}
\vspace{2cm}

\begin{center}
\larger[2]\bf
The Penrose transform on conformally Bochner-K\"ahler manifolds
\end{center}

\vspace{1cm}

\begin{center}
	Yoshinari Inoue
\end{center}

\vspace{5mm}

\begin{center}
{\smaller
\noindent
Department of Mathematics, Graduate School of Science,\\
 Kyoto University, Kyoto 606--01, Japan\\
\ifMA
 Fax: +81--75--753--3711\\
 E--mail: {\tt inoue@kusm.kyoto-u.ac.jp}
\else
 E--mail: {\tt inoue@kusm.kyoto-u.ac.jp}\\
 \phantom{F0}
\fi
}
\end{center}

\vspace{-5mm}

\begin{abstract}
We give a generalization of the Penrose transform on Hermitian
manifolds with metrics locally conformally equivalent to
Bochner-K\"ahler metrics. We also give an explicit formula for
 the inverse transform.
\end{abstract}

\maketitle

{\smaller[2]\noindent
 \subjclassname: 32L25 (Primary), 53C55, 32C35 (Secondary)}

\vspace{2mm}

%
\section*{Introduction}

The twistor space of a four-dimensional Riemannian manifold
was first introduced by Atiyah, Hitchin and Singer in \cite{[A.H.S]}.
The Penrose transform on it was given by Hitchin in \cite{[H]}.
The inverse Penrose transform was given explicitly by
Woodhouse in \cite{[W]}.

In \cite{[O.R]}, O'Brian and Rawnsley extended
the notion of twistor space to other manifolds such as
even dimensional Riemannian manifolds, Hermitian manifolds,
Quaternionic K\"ahler manifolds, etc.
In \cite{[M]}, Murray gave the Penrose transform on twistor spaces
of even dimensional Riemannian manifolds.
An explicit formula for the inverse Penrose transform
in this case was given by the author in \cite{[I3]}.

As for twistor spaces of Hermitian manifolds,
the author gave the Penrose transform
and its inverse transform  on the twistor spaces
of $\C^n$ in his unpublished work  \cite{[I2]}.
They related a certain cohomology group on a twistor space
to the space of harmonic forms of the Dolbeault
complex.
An explicit formula of the inverse transform was also given.

In the present paper, we
extend the result of \cite{[I2]} to
Hermitian manifolds whose metrics are locally conformally
equivalent to
Bochner-K\"ahler metrics, which we call
{\em conformally Bochner-K\"ahler manifolds}.
Since the Penrose transform involves
cohomology groups of positive degree, the integrability condition for the
almost complex structure of the twistor space is inevitable.
Hence it is a complete generalization of \cite{[I2]}
except low dimensional examples of non-integrable almost Hermitian
manifolds, such as $S^6$ with a well-known almost complex structure.
Since our assumption on the metric of the base manifold
 is slightly different
from the one  in \cite{[O.R]},
there might be
new exceptional examples
with real dimension
eight (Theorem \ref{thm:int}).

Since we have no reason to exclude the conformally flat case,
all results and proofs in \cite{[I2]}  are
contained in this paper.

Let $X$ be a conformally Bochner-K\"ahler manifold of dimension $n>2$.
Let $k$ be an integer between $0$ and $n$.
Then the twistor space $Z_k(X)$
is a complex manifold (\cite{[O.R]}, \cite{[I4]}).
Let $H$ be the hyperplane  bundle as a Grassmannian bundle.
Let $h$ be a non-negative integer and $V$ a vector bundle
with connection
on $X$ such that $H^{-n-h}\otimes p^*V$ has a $(1,0)$-connection
where $p: Z_k(X) \rightarrow X$ is the projection map.
The trivial bundle does not satisfies this condition
if the metric of $X$ is not locally conformally equivalent
to a flat metric,
which is quite different in the Riemannian case.
The Penrose transform is a one-to-one correspondence
between
the cohomology group $H^{k(n-k)}(Z_k(X), \cO(H^{-n-h}\otimes p^*V))$
and the solution space of the Laplacian $\mathcal{D}_0$ defined
in Definition~\ref{defn:D0} if $h=0$,
or the space of harmonic forms of the Dolbeault complex:
\begin{align*}
\begin{split}
	\Gamma(\Wedge^{0,k-1}X\otimes \hS^{h-1} \Wedge^{0,k}X\otimes
	 K_X^{-k}\otimes V) & \\
	  \buildrel \db \over
	 \longrightarrow
	\Gamma(\hS^h \Wedge^{0,k}X\otimes &K_X^{-k}\otimes V)
	  \buildrel \db \over \longrightarrow\\
	 &\Gamma(\Wedge^{0,k+1}X\otimes \hS^{h-1} \Wedge^{0,k}X\otimes
	  K_X^{-k}\otimes V),
\end{split}
\end{align*}
if $h>0$,
where $\hS^{h} \Wedge^{0,k}$ is the principal component of
$S^h\Wedge^{0,k}$ defined in Definition~\ref{defn:prod}
(Theorem \ref{z2x} and Theorem \ref{thm_x2z}).

The explicit formula for the inverse Penrose transform
given in Definition \ref{defn:inv}
has quite similar form as in the Riemannian case (\cite{[I3]}).
Namely put
$$
	F(x) = \sum_{i=0}^{\infty} \frac{x^i}{(i!)^2},
 $$
and
\begin{align*}
 	f(x) &= (k+h-1)! F^{(k+h-1)}(x),\\
	g(x) &= (n-k+h-1)! F^{(n-k+h-1)}(x).
\end{align*}
Let $j$, $D_\alpha$, $D_\beta$ and $D_\gamma$ be operators
defined in \S\ref{sec:inv}.
Then the inverse Penrose transform is given by:
$$
\begin{matrix}
	\AAA:& \Gamma(X, \hS^h\Wedge^{0,k}X\otimes
		K_X^{-k}\otimes V) &\longrightarrow
		&\Omega^{0,k(n-k)}(Z_k(X), H^{-n-h}\otimes p^*V)\\[3mm]
		& \phi &\longmapsto
		& f(D_\beta + \Dr) g(D_\alpha)
			\jjj
\end{matrix}
 $$
where the expression is well-defined because $\Da$ $\Db$ and $\Dr$
are even operators and $\Db$ and $\Dr$ are mutually commutative.

Let us explain briefly the contents of this paper. We define in
 \S\ref{sec:twh} the twistor spaces of an almost Hermitian manifold
and study their integrability conditions.
In \S\ref{sec:curv} we give the field equations
on conformally
Bochner-K\"ahler manifolds which appear in the Penrose transform.
We also give a few formulas about the curvature tensor
which are used in \S\ref{sec:inv}.
In \S\ref{s:pt} we define the Penrose transform on conformally
Bochner-K\"ahler
manifolds. In \S\ref{sec:inv} we prove the
surjectivity of the Penrose transform by constructing explicitly
 the inverse correspondence.
We summarize in Appendix the operators which are
defined during the computation in \S\ref{sec:inv}.

%
\section{The twistor spaces of an almost Hermitian manifold
and their integrability conditions
}\label{sec:twh}

Let $X$ be  an $n$-dimensional
almost Hermitian manifold, that is, $X$ is a $2n$-dimensional
Riemannian manifold with
a compatible almost complex structure.
In this section we define twistor spaces
 $Z_k(X), k=0,\dots,n$ of $X$ as submanifolds of the
Riemannian twistor space $Z(X)$. It is
modification of the definition by O'Brian and Rawnsley in \cite{[O.R]}.

	Let
$$
	\D = \langle \theta_I \mid I < (1,\dots,n) \rangle_\C
 $$
be a spin module,
where $I<(1,\dots,n)$ means that $I$ is a subsequence of $(1,\dots,n)$.
We introduce similar notation for multi-indices as in \cite{[I3]}.
We regard  multi-indices $I,J,\dots$ as finite sequences
of possibly duplicate elements of $\{1,\dots,n\}$, and denote by
$IJ$ the composition of sequences $I$ and $J$.
Let us define relations among ${\theta_I}$'s as:
\begin{align*}
	  \theta_{I_1aaI_2} &= - \theta_{I_1I_2},\\
	  \theta_{I_1abI_2} &= - \theta_{I_1baI_2}, \quad\hbox{for $a\ne b$.}
\end{align*}
Then, for any multi-index $I$,  there is a unique subsequence
$I_0$ of $(1,\dots,n)$ such that
$\theta_I = \theta_{I_0}$ or $\theta_I = - \theta_{I_0}$.
If a multi-index $I$ is regarded as a set,
 it is meant as the set of numbers
contained in $I_0$.
Let $|I|$ denote the length of $I_0$.

With this notation, we define an action of
$\R^{2n} = \langle f_i \mid i=1,\dots,2n \rangle_\R$ on $\D$
 as follows. For $a=1,\dots,n$,
\begin{align*}
	  f_a \theta_I &= \theta_{aI},\\
	  f_{n+a} \theta_I &=
		\begin{cases}
		  \im \theta_{aI},	&\text{if $a \notin I$},\\
                  -\im\theta_{aI}, &\text{if $a \in I$.}
		\end{cases}
\end{align*}
This action is canonically extended to a $\Clif(\R^{2n})$-action on $\D$.
Since $\spin(2n)$ is a subspace of $\Clif(\R^{2n})$, we have
a $\spin(2n)$-action on $\D$. (See \cite{[I1]} for details.)

  If we identify $\R^{2n}$ with $\Cn$ by the complex structure
defined by
$$
	J = \begin{pmatrix}
		0 & -I_n\\ I_n & 0
	     \end{pmatrix},
 $$
then $\R^{2n}\otimes \C$ is decomposed into two submodules:
\begin{align*}
	\C^n \simeq
	\langle e_a =  \frac{1}{2} (f_a - \im f_{n+a})
		\mid a=1,\dots,n \rangle_\C, \quad
	{\ol \C}^n \simeq
	\langle e_{\ol a} = \overline{e_a}
		\mid a=1,\dots,n \rangle_\C.
\end{align*}
Actions of these vectors on $\D$ are computed as:
\begin{equation}\label{eq:cm}
	e_a\theta_I =
	  \begin{cases}
		\theta_{aI}, &\text{if $a\notin I$}\\
		0, &\text{if $a\in I$}
	  \end{cases} \qquad
	e_{\ol a}\theta_I =
	  \begin{cases}
		0, &\text{if $a\notin I$}\\
		\theta_{aI}, &\text{if $a\in I$}
	  \end{cases}
\end{equation}
If we consider $\un(n)$ as a subspace of $\Clif(\R^{2n})$, we have
$$
	\un(n) \subset \langle e_a\cdot e_{\ol b}
	\mid a,b = 1,\dots, n \rangle_\C.
 $$
Hence the irreducible decomposition of $\D$ as a $\un(n)$-module
is given by:
\begin{align*}
	  \D = \Oplus_{k=0}^n \D^k,\qquad
	  \D^k = \langle \theta_I \mid |I| = k \rangle_\C,
\end{align*}
where we have
\begin{equation}\label{eq:rep_delta_k}
	\D^k \simeq (\Wedge^k \Cn) \otimes (\Wedge^n \Cn)^{-1/2}.
\end{equation}

We fix an integer $k$ such that $0 \le k \le n$.
Let $Z$ denote the $\Pin(2n)$-orbit of $[\theta_\emptyset]\in \Proj(\D)$,
which is a complex submanifold of $\Proj(\D)$.
Now we study the submanifold
\begin{align*}
	Z_k = Z \cap \Proj(\D^k).
\end{align*}
We see that $Z_0$ and $Z_n$ are one point corresponding
to the original complex structure and its complex conjugate
structure, respectively.
Henceforth we assume $0<k<n$.

To simplify the expression, we stretch the rules of notation
of multi-indices as follows.
Let $(Z^I)_{I<(1,\dots,n)}$ be a coordinate system
corresponding to $(\theta_I)_I$. Then, by regarding it as a cospinor,
we can define $Z^{aI}$ in a similar way.
Since we want to consider $Z^{aI}$ as a homogeneous coordinate
of $\Proj(\D^k)$, the length of the multi-index $aI$ is significant.
Therefore we write $aI$ (resp. $\ol aI$) if $\abs{aI}= \abs{I}+1$
 (resp. $\abs{aI} = \abs{I} -1$).
In the same way, we write $\ol a\ol I$ (resp. $a\ol I$)
if $\abs{aI} = \abs{I} +1$ (resp. $\abs{I}-1$).
These rules are
 inspired by the formula \eqref{eq:cm} for Clifford multiplication.
Owing to it, as we will see below, we can
omit the range of summation in most cases.

  Let $J$, $K$ be multi-indices of length $|J|=k+1$, $|K|=k-1$.
By \cite{[I3]} Lemma~2.6~(1), we have
$$
	\sum_{a\in J\setminus K} Z^{\ol aJ}Z^{aK} +
	\sum_{a\in K\setminus J} Z^{aJ}Z^{\ol aK}
	= 0,
 $$
on $Z$. Since the second term vanishes on $\Proj(\D^k)$, we have
\begin{align}\label{eq:plucker}
	Z^{\ol aJ}Z^{aK} = 0,
\end{align}
on $Z_k$, which is nothing but the Pl\"ucker relation.
Therefore $Z_k$ coincides with the
Grassmannian manifold $G_{k,n}$.

 Hence we consider a diagram:
\begin{equation}\label{zk_in_z}
\begin{CD}
	  G_{k,n} @>>> \Proj(\D^k)\\
	  @VVV @VVV \\
	  Z @>>> \Proj(\D)
\end{CD}
\end{equation}
where upper (lower) horizontal array is a $\Un(n)$- ($\Or(2n)$-)
equivariant map.

	Let $X$ be an $n$-dimensional almost Hermitian manifold.
Let $P$ denote the principal
bundle of $X$. Then we define the $k$th twistor space as:
\begin{align*}
	  Z_k(X) = P \times_{\Un(n)} G_{k,n}.
\end{align*}
Let $p$ denote the projection map:
$$
	p: Z_k(X) \longrightarrow X.
 $$

If $X$ has a spin structure,
we can define the spin-hyperplane bundle $H_0$ over $Z_k(X)$ by pulling back
the hyperplane bundle over $Z(X)$ (or equivalently $\Proj(\D^k(X))$).
On the other hand, if we regard $Z_k(X)$
as a Grassmannian bundle $G_k(T^{(1,0)}X)$, it is natural to consider the
Grassmannian-hyperplane bundle $H$,
which can be defined even if $X$ does not have a  spin structure.
Hence we use it to study the Penrose transform
and call it the hyperplane bundle of $Z_k(X)$.
By \eqref{eq:rep_delta_k}, two line bundles are related by
the isomorphism:
$$
	H \simeq H_0 \otimes p^*K_X^{\frac 12},
 $$
where $K_X$ is the canonical bundle of $X$.

By \eqref{zk_in_z}, $Z_k(X)$ is a subbundle of $Z(X)$. Since
$G_{k,n}$ is a complex submanifold of $Z$, we can define an almost
complex structure of $Z_k(X)$  similarly as
in the case of $Z(X)$.

We take a local orthonormal frame $(e_a)_a$ of $T^{(1,0)}X$.
Let $Z^J$ denote
 a corresponding Pl\"ucker coordinate of $Z_k(X)$.
Then we can define local fiber coordinates of $Z_k(X)$ as
follows.
Let $I$ be a multi-index of length $k$. Then
$$
	w_{ij}=Z^{\ol\imath jI}/Z^I,\qquad i\in I,\quad j\not\in I
 $$
are local  fiber coordinates on
$$
	U_I=\{(Z^J)_J \in G_{k,n} \mid Z^I\ne 0\}.
 $$
We denote  the center of the coordinate chart $U_I$ as $z_0$.
For a general multi-index $J$ of length $k$,
 let $z^J$ be the function on $U_I$ defined as
$z^J = Z^J/Z^I$.
Let $\omega^J_K$ be the connection form
of $\Wedge^kT^{(1,0)}X$ induced by the Levi-Civita connection.
Let $(e^a)_a$ be the dual frame of
$\Wedge^{1,0}X$.

By \cite{[I3]} Lemma 2.9,
 the $(1,0)$-forms of $Z_k(X)$ can be explicitly defined as follows.

\begin{defn}\label{def:acs}
The total space of ${H^*}^\times\subset \Wedge^kT^{(1,0)}X$ has
an almost complex structure
whose space of $(1,0)$-forms
is spanned by
\begin{align*}
	Z^{aJ} e^{\ol a},& \qquad |J|=k-1,\\
	Z^{\ol aJ} e^a,& \qquad |J|=k+1,\\
	dZ^J + \omega^J_K Z^K,& \qquad \abs{J}= k.
\end{align*}
This induces an almost complex
structure on $Z_k(X)$ whose space of $(1,0)$-forms is spanned by
\begin{align*}
	z^{aJ} e^{\ol a},& \qquad |J|=k-1,\\
	z^{\ol aJ} e^a,& \qquad |J|=k+1,\\
	\widehat {dw_{ij}},& \qquad i \in I,\quad j \not\in I,
\end{align*}
where
$$
	\widehat {dw_{ij}} = dw_{ij} + z^J\omega^{\ol \imath jI}_J
	 - w_{ij} z^J \omega^I_J.
 $$
\end{defn}

\begin{rem}
The integrability of $Z_k(X)$
does not implies that of ${H^*}^{\times}$  in general.
Hence we should twist the hyperplane bundle by the pull-back of
a non-trivial vector bundle on $X$ satisfying a certain
condition on the curvature.
The condition shall be calculated in the next section.
\end{rem}

Since the Penrose transform in the two-dimensional (i.e. real dimension
four) case  is already given by Hitchin in \cite{[H]}, we consider
henceforth the case $n>2$.

O'Brian and Rawnsley
considered
in \cite{[O.R]}
 an almost complex structure
defined as above but use
 a general connection
which makes the almost complex structure of $X$ and the
metric tensor of $X$ covariant
constant.
To avoid unnecessary generality, we restrict ourselves to
consider the Levi-Civita connection.
Then the extra conditions on the connection
does not seem very good, since it
 is not conformally invariant.
We want to avoid this by considering instead the condition that
$Z_k(X)$ is an almost complex submanifold of $Z(X)$.

First, we give a definition of Bochner-K\"ahler manifolds.
Let $X$ be a K\"ahler manifold.
Let $(e_a)_a$ be a local orthonormal
frame of the holomorphic tangent bundle $T^{(1,0)}X$.
Let $(e^a)_a$ and $(e^{\ol a})_{\ol a}$ be the dual frame and
its complex conjugate frame, respectively.
Let
$$
	R^b_a e_b =
	  R_{c\ol da}^b e_b\otimes e^c\wedge e^{\ol d}
 $$
 be the curvature tensor of $T^{(1,0)}X$
induced by the Levi-Civita connection.
Put
$$
	\Curv abcd = \frac 12 R_{a\ol bc}^d.
 $$
The Ricci curvature and the scalar curvature are defined as
\begin{align*}
	\Ric{a}{b} &= \frac{2}{n+2} \Curv abcc,\\
	r &= \frac{1}{n+1} \Ric aa.
\end{align*}
The multiplicative constants in the definitions are just for simplicity.

The curvature tensor of a K\"ahler manifold
 has three irreducible components:
the scalar curvature, the traceless Ricci curvature and the Bochner
curvature.
Hence if the Bochner tensor of a manifold  vanishes, its
curvature tensor
can be written explicitly by the scalar and Ricci curvatures.
Therefore, we can define a Bochner-K\"ahler manifold by
that explicit formula as in \cite{[Ka]}.
\begin{defn}\label{defn:bk}
A K\"ahler manifold $X$ is called {\em a Bochner-K\"ahler manifold\/} if
its curvature tensor is written as:
$$
	\Curv abcd = \frac 12 (\Ric ab \delta_c^d
			+ \Ric ad \delta_c^b
			+ \Ric cb \delta_a^d
			+ \Ric cd \delta_a^b)
		- \frac 12 r (
			  \delta_a^b\delta_c^d
			+ \delta_a^d\delta_c^b).
 $$
\end{defn}

Now we return to consider a general almost Hermitian manifold $X$.
We say that an almost Hermitian
manifold is {\em conformally Bochner-K\"ahler\/}
if it is a complex manifold and
its metric is locally conformally equivalent to a Bochner-K\"ahler
metric.

\begin{defn}
The pair $(n,k)$ is called {\em exceptional\/} if it is equal to
$(3,1)$, $(3,2)$ or $(4,2)$.
\end{defn}
Essential part of
the following theorem was  given by O'Brian and Rawnsley
in \cite{[O.R]}.

\begin{thm}\label{thm:int}
Let $X$ be an $n$-dimensional almost Hermitian manifold with $n>2$.
Let $k$ be an integer between $0$ and $n$.
Suppose that $Z_k(X)$ is a complex manifold and also
 an almost complex
submanifold of $Z(X)$.
\begin{enumerate}
\item
If $(n,k)$ is not exceptional,
then $X$ is a conformally Bochner-K\"ahler manifold.
\item
  Suppose $(n,k)$ is exceptional
and $X$ is a complex manifold.
Then $X$ is a conformally Bochner-K\"ahler manifold.
\end{enumerate}
\end{thm}

\begin{proof}
We give the condition at $z_0$ explicitly.

  The condition that $Z_k(X)$ is an almost complex submanifold
of $Z(X)$ is
\begin{align*}
	\omega^{\ol a\ol bI}_I, \omega^{abI}_I \in \Wedge^{1,0}_{z_0},
\end{align*}
where $\omega^J_K$ is meant as the connection form of $\D(X)$
and $\Wedge^{1,0}_{z_0}$ is the space of $(1,0)$-forms at $z_0$.
Let
$$
	  T^{\ol a}_{c,b} e^c + T_{\ol c,b}^{\ol a} e^{\ol c} =
	  \langle e^{\ol a}, \omega(e_b)\rangle
 $$
be torsion tensors.
Note that this is {\em not\/} the torsion as a connection of $TX$,
which vanishes because we consider the Levi-Civita connection.
Then the condition is rewritten as:
\begin{align}
	T_{c,b}^{\ol a} = 0,&\qquad a,b,c\in I,\quad a\not=b,
	  \label{eq:t1}\\
	T_{\ol c,b}^{\ol a} = 0,
		&\qquad a,b\in I,\quad c\not\in I,\quad a\not=b,
	  \label{eq:t2}\\
	T_{\ol c,b}^{\ol a} = 0,
		& \qquad a,b\not\in I,\quad c\in I,\quad a\not=b,
	  \label{eq:t3}\\
	T_{c,b}^{\ol a} = 0,& \qquad a,b,c\not\in I,\quad a\not=b.
	  \label{eq:t4}
\end{align}

By \eqref{eq:t2} and \eqref{eq:t3},
$T_{\ol c,b}^{\ol a} = 0$ for distinct $a$, $b$ and $c$.
This means that there is a tensor $(t_a)_a$
such that
\begin{align*}
	T_{\ol c,b}^{\ol a} = \delta^c_at_b - \delta^c_bt_a.
\end{align*}

By \eqref{eq:t1} and \eqref{eq:t4},
we have
$$
	T_{c,b}^{\ol a} = 0,
 $$
if $k\ge3$ or $n-k\ge3$,
which is satisfied if and only if $(n,k)$ is not exceptional.
This tensor is nothing but the Nijenhuis tensor. Hence its vanishing
implies precisely that $X$ is a complex manifold.

  Therefore, under the assumption of the theorem,
for the K\"ahler form $\Omega$, we have
$$
	d\Omega = 4 t\wedge\Omega,
 $$
where $t\in\Gamma(T^*X)$ is the {\em real\/} section corresponding
to $t_a e^a$.
Then $dd\Omega = 0$ implies that $t$ is a closed form
and there exists a locally defined
scalar function $s$
such that $ds = t$.
Hence the Hermitian metric obtained from the original metric
by multiplying $\exp(-4s)$ is torsion free.

 This implies that
$X$ is a complex manifold and its metric is locally
conformally equivalent to a K\"ahler metric.

Hence we can assume that $X$ is a K\"ahler manifold and
the connection form of $T^{(1,0)}X$
vanishes at a point $x_0\in X$.
Then we have
\begin{align}\label{eq:ddwa}
	\left. d \dwa \right|_{(x_0,z_0)} = R^j_i.
\end{align}
Hence the integrability condition of the almost complex structure
of $Z_k(X)$ at $(x_0,z_0)$ is
$$
	R_{k,\ol l, i}^j = 0,\qquad
	  i, k \in I,\quad j,l \not\in I.
 $$
This implies precisely
the vanishing of the Bochner tensor by Definition~\ref{defn:bk}.
\end{proof}

\begin{expl}
In the exceptional case,
$S^6$ with
the almost complex structure
induced by the Cayley algebra
has the integrable twistor spaces  $Z_k(S^6)$, $k=1,2$,
which were described in \cite{[O.R]}.
The twistor space $Z_+(S^6)$ is identified with the
homogeneous space $\SO(7)/\Un(3)$.
The twistor space $Z_2(S^6)$ is identified with $G_2/\Un(2)$
where $G_2$ is the subgroup of $\Spin(7)$ which
consists with transformations preserving the almost complex
structure of $S^6$.

We give here an explicit description of $Z_2(S^6)$
as a submanifold of $Z_+(S^6)$.
The space $Z_+(S^6)$ is identified with the six-dimensional
complex hyperquadric:
$$
	Q_6 = \{ [Z^I]_{I<(1,2,3,4)}\mid
	  Z^\emptyset Z^{1234}-Z^{12}Z^{34}+Z^{13}Z^{24}-Z^{14}Z^{23} = 0 \}
 $$
(see \cite{[I1]} for detail).
The defining equation
of $Z_2(S^6)\subset Q_6$
obviously depends on the choice of the almost complex structure
of $S^6$, or equivalently the embedding $G_2\subset \Spin(7)$.
For example, if $G_2$ is the isotropic subgroup at the
(co)spinor $\theta^\emptyset + \theta^{1234}$,
$Z_2(S^6)$ can be written  as
$$
	Z_2(S^6) = \{ [Z^I]\in Q_6 \mid Z^\emptyset + Z^{1234} = 0\}.
 $$
Hence it is a five-dimensional complex hyperquadric.
\end{expl}

\begin{rem}
This theorem suggests that there might exist
a four-dimensional non-integrable
almost Hermitian manifold $X$ having the integrable second
twistor space $Z_2(X)$.
The result of this paper is expected to extend
to these non-integrable almost Hermitian manifolds.
\end{rem}

  The almost complex structure of $Z_k(X)$ can also be defined
by using the distribution of some first order differential
operator.
This gives immediate correspondences between
$0$-th cohomology groups and some
field equations as follows.

\begin{defn}\label{defn:prod}
Let $E_1$ and $E_2$ be irreducible $\Un(n)$-modules.
For $i=1,2$,
let $v_i\in E_i$ be a highest weight vector with weight $\lambda_i$.
Then we write
$E_1\hotimes E_2$ as
a unique irreducible submodule of $E_1\otimes E_2$
having a highest weight vector $v_1\otimes v_2$
with weight $\lambda_1 + \lambda_2$.
 Similarly we write $\hS^h E_1$
as a unique irreducible
submodule of $S^h E_1$ with highest weight $h\lambda_1$.
\end{defn}

Put $E = \Wedge^{k,0}\hotimes \Wedge^{1,0}\oplus
\Wedge^{k,0}\hotimes \Wedge^{0,1}$.
Then, by the Littlewood-Richardson rule
for the irreducible decomposition of the tensor representation,  we have
$$
	\Wedge^{k,0}\otimes (\Wedge^{1,0}\oplus \Wedge^{0,1})
	= \Wedge^{k-1,0}\oplus \Wedge^{k+1,0} \oplus E.
  $$
By composition of the covariant derivative and  the projection,
we define a first order differential operator:
\begin{align*}
	\bar{\mathcal D} : \Gamma(\Wedge^{k,0}X) \longrightarrow
	  \Gamma(E(X)).
\end{align*}
Then the distribution $V(\bar{\mathcal D})$
 of $T\Proj(\Wedge^{k}T^{(1,0)}(X))\otimes \C$
defined in \cite{[A.H.S]} and \cite{[I1]}
has minimum rank on $Z_k(X)$ and gives on it the almost complex
structure of Definition~\ref{def:acs}.

 Let $V$ be a vector bundle on $X$ with connection.
Let $h$ be a non-negative integer.
  As in \cite{[H]} \S2  and \cite{[I1]} \S9, we have an immediate
correspondence between $H^0(Z_k(X), \cO(H^h\otimes V))$ and the solution
space
of the  equation $\bar{\mathcal D}_h\phi = 0$ where
$\bar{\mathcal D}_h$ is the differential operator
induced by the covariant derivative
similarly as $\bar{\mathcal D}= \bar{\mathcal D}_1$.
\begin{align*}
	\bar{\mathcal D}_h : \Gamma(\hS^h\Wedge^{k,0}X \otimes V)
	 \longrightarrow
	\Gamma ( (\Wedge^{1,0}X\hotimes\hS^h\Wedge^{k,0}X\oplus
	  \Wedge^{0,1}X\hotimes\hS^h\Wedge^{k,0}X)
	  \otimes V).
\end{align*}

Note that the correspondence is valid even if $Z_k(X)$ is not
a complex manifold.

%
\section{Field equations on a conformally
	Bochner-K\"ahler manifold}\label{sec:curv}

Let $n$ be an integer greater than two.
Let $X$ be an $n$-dimensional conformally Bochner-K\"ahler
 manifold.
We fix an integer $k$ such that $0 < k < n$.
In this section, we introduce
 field equations on $X$ which appear in the Penrose transform.

Since the construction of the field equations are local,
we can assume that the metric of $X$
is a Bochner-K\"ahler metric.
We note that
the K\"ahler condition determines a metric in the conformal class
up to a locally constant scalar factor.
Therefore the situation is much simpler than that in the
case of Riemannian manifolds (\cite{[I3]}).

Our main concern  is a relation between the solution space
of a field equation
and a cohomology group
of positive degree
 with coefficients in the vector bundle
$H^{-n-h}\otimes p^*V$,
where $h$ is a non-negative integer
and $V$ is a vector bundle on $X$ with connection.
Hence the vector bundle $H^{-n-h}\otimes p^*V$
 should have a natural holomorphic structure.

Now we give an explicit condition on $V$ and
compute the curvature tensor of
$\hS^h\Wedge^{0,k}X\otimes K_X^{-k}\otimes V$,
on which the field equation is defined.

Let $R_a^b$ be the curvature tensor of the holomorphic tangent bundle of $X$.
Then
the curvature tensor of $\Wedge^{0,k}X$ is
\begin{align*}
	(R\psi)_{\ol I}
	  = -R_a^b \psi_{\ol ab\ol I}.
 \end{align*}
By the definition of the Ricci tensor,
the curvature tensor of the canonical bundle $K_X$ is
\begin{equation}\label{curv_k}
	- (n+2) \Ric ab e^a \wedge e^{\ol b}.
\end{equation}
Now we compute the condition on $V$. Let $R^{I}_{J}$ be
the curvature tensor of $\Wedge^{k}T^{(1,0)}X$. Then, by
Definition \ref{defn:bk},
we have
\begin{align*}
	 R^{I}_{I} &= \sum_{a\in I} R^a_a
	= (k+1) \Ric ab e^a\wedge e^{\ol b}
	 + \text{$(1,1)$-forms at $z_0$}.
 \end{align*}
If a two-form is of type $(1,1)$ for all almost
complex structures corresponding to points of $G_{k,n}$,
then it
should be a scalar multiple of the K\"ahler form.
Here we use the assumption $n>2$.
Hence we have proved the following lemma.
\begin{lem}\label{lem:curv_v}
Let $V'$ be a vector bundle on $X$ with connection.
Then,
for an integer $l$,
the vector bundle $H^l\otimes p^*V'$ has a natural holomorphic
structure if and only if the curvature of $V'$ is written as:
$$
	l(k+1)\Ric ab e^a\wedge e^{\ol b} + e^a\wedge e^{\ol a} \El,
 $$
where $\El$ is an endomorphism of $V'$.
\end{lem}
\begin{rem}
By \eqref{curv_k}, we can use
$H^l\otimes p^*K_X^{-l\frac{k+1}{n+2}}$
as  a locally defined holomorphic line bundle.
\end{rem}

\begin{expl}
We consider the condition
that $H^{-n-h}\otimes p^*V$ is a holomorphic line bundle
when $X=\Hypb^l\times\Proj^{n-l}$ and
$V=L(a,b) = p_1^*(H^a)\otimes p_2^*(H^b)$, where $p_i$ is the projection
to the $i$th component.
We have
$$
	K_X \simeq L(-(l+1),-(n-l+1)).
 $$
Since the curvature tensor of $L(1,-1)$
is a scalar multiple of the K\"ahler form of $X$,
by \eqref{curv_k},
the condition for $V$ is
$$
	a + b =  -(n+h)(k+1).
 $$
Hence we have
$$
	K_X^{-k}\otimes V = L(a',b'),\qquad
		a' + b' = -  h(k+1)  - (n-2k).
 $$
\end{expl}

Suppose that the vector bundle $V$ has the curvature such that
$H^{-n-h}\otimes p^*V$ has a holomorphic structure.
Then, for $\phi\in\Gamma(\hS^h\Wedge^{0,k}\otimes K_X^{-k}\otimes V)$,
the curvature tensor is written as:
\begin{equation}\label{eq:spin_curv}
 \begin{split}
	(R\phi)_\iii
	   &= - \sum_{j=1}^h R_b^a\phi_{\ol{I}_1,\dots
		,\ol ba\ol{I}_j,\dots,\ol{I}_h}
		- ( h(k+1) + n - 2k)\Ric ab e^a\wedge e^{\ol b}
	  \phi_\iii\\
	&\quad
	  +  e^a \wedge e^{\ol a}(\El\phi)_\iii,
  \end{split}
 \end{equation}
where $\El$ is an endomorphism of $V$.
\newcommand{\dirac}{{\mathcal D}}

We assume for a while that $h\ge1$.
Let $\dirac_h$ denote the Dirac operator.
We obtain its harmonic section from a harmonic form of
$\db+ \db^*$ by the following commutative diagram:
\begin{align*}
\begin{CD}
	\Gamma(\hS^h\Wedge^{0,k}X\otimes K_X^{-k}\otimes V) @>\db+\db^*>>
	 \Gamma ((\Wedge^{0,k-1}\oplus \Wedge^{0,k+1})\otimes
	    \hS^{h-1}\Wedge^{0,k}X\otimes K_X^{-k}\otimes V)\\
	@VVV  @VVV\\
	\Gamma(S^h \D^{\sigma(k)}X\otimes K_X^{-\frac h2-k}\otimes V)
		 @>\dirac_h>>
	  \Gamma(\D^{\sigma(k+1)}\otimes S^{h-1}\D^{\sigma(k)}X\otimes
		 K_X^{-\frac h2-k}\otimes V)
\end{CD}
\end{align*}
where $\sigma(k)$ is $+$ (resp. $-$) if $k$ is even (resp. odd).
This diagram also shows that $\db^*$ can be defined as a global
operator if we restrict the metric locally chosen
to be K\"ahler.
Henceforth we use the operator $\db^*$ in this sense.

For a positive integer $a$, we consider
$$
	\nabla^a: \Gamma(V) \longrightarrow \Gamma(V\otimes
			\Otimes^a T^*X),
 $$
where the connection of $T^*X$ is the one induced from
the Levi-Civita connection.
Evaluation by tangent vectors is defined inductively as:
$$
	\langle \nabla^a\phi, e_a\dots e_1\rangle =
	  \langle \nabla_{e_1}
	    \nabla^{a-1}\phi, e_{a}\dots e_2\rangle.
 $$
Let $\phi$ be a harmonic form of
$\hS^h\Wedge^{0,k}X\otimes K_X^{-k}\otimes V$.
Then by considering it as a harmonic section of
$S^h\D^{\sigma(k)}X\otimes K_X^{-\frac h2 -k}\otimes V$,
we have
\begin{align*}
	\frac 14 (\dirac_h\dirac_h\phi)_{\ol{I},\dots,\ol{I}}
	&=
	  \langle
	    \nabla^2 \phi, e_a e_{\ol b}
	  \rangle_{\ol a b\ol{I},\dots,\ol{I}}
	+
	  \langle
	    \nabla^2 \phi, e_{\ol b} e_a
	  \rangle_{b \ol a \ol{I},\dots,\ol{I}}\\
	&=
	- \langle
	    \nabla^2\phi, e_a e_{\ol a}
	  \rangle_{\ol{I},\dots,\ol{I}}
	+ \langle
	    R\phi, e_{\ol b} e_a
	  \rangle_{b \ol a \ol{I},\dots,\ol{I}}\\
	&= 0.
 \end{align*}
Hence, by \eqref{eq:spin_curv}, we have
\begin{align*}
	  \langle
	    \nabla^2\phi,
	      e_a  e_{\ol a}
	  \rangle_{\ol{I},\dots,\ol{I}}
	&= - 2 \Curv abcd \phi_{\ol cd b\ol a\ol{I},\dots,\ol{I}}
	   - 2(h-1) \Curv abcd
	     \phi_{b\ol a\ol{I},\ol cd \ol{I},\dots,\ol{I}}\\
	&\quad - (h(k+1)+n-2k) (S_0\phi)_{\ol{I},\dots,\ol{I}}
		- (n-k) (\El\phi)_{\ol{I},\dots,\ol{I}}.
\end{align*}
Here $S_0$ is the endomorphism
\begin{equation}\label{eq:S0}
	S_0\phi
	  = \Ric ab \phi_{b\ol a\ol{I}_1,\dots,\ol{I}_h} e^\iii,
\end{equation}
where $e^\iii$ denotes the image of
$e^{\ol{I}_1}\otimes\dots \otimes e^{\ol{I}_h}$ by
the projection $\Otimes^h \Wedge^{0,k}X\rightarrow \hS^h\Wedge^{0,k}X$.
By using Definition \ref{defn:bk}, we compute
\begin{gather*}
	2 \Curv abcd \phi_{\ol cd b \ol a\ol{I},\dots,\ol{I}} = 0,\\
	2 \Curv abcd \phi_{b\ol a\ol{I},\ol cd\ol{I},\dots,\ol{I}}
	 = (n-2k)(S_0\phi)_{\ol{I},\dots,\ol{I}}
	  +(n-k)(n-k+1)r\phi_{\ol{I},\dots,\ol{I}}.
\end{gather*}
Hence, by using the irreducibility of $\hS^h\Wedge^{0,k}$, we have
\begin{equation}\label{eq:lap}
\begin{split}
	  \langle
	    \nabla^2\phi,
	      e_a  e_{\ol a}
	  \rangle
	&= - h(n-k+1) S_0\phi
	 - (h-1)(n-k)(n-k+1) r \phi\\
	&\quad
	 - (n-k)\El\phi.
\end{split}
\end{equation}

In order to define the field equation when $h=0$,
we use a similar method which was used
to define the conformally invariant
Laplacian by Hitchin in \cite{[H]}.

Let $V$ be a vector bundle with connection
on $X$ such that $H^{-n}\otimes p^*V$
has a holomorphic structure.
Let $V'$ be the locally defined line bundle $K_X^{-\frac{k+1}{n+2}}$.
Then, by Lemma \ref{lem:curv_v},
 $H\otimes p^*V'$ is a holomorphic line bundle.
Let $\psi\in\Gamma(\Wedge^{k,0}X\otimes V')$ be a section
satisfying the equation $\bar{\mathcal D}\psi = 0$.
Then we have
\begin{align*}
	\langle \nabla^2\psi, e_a e_{\ol a} \rangle_I
	  &= \sum_{a\not\in I}
		\langle R\psi, e_a e_{\ol a} \rangle_I\\
	  &= (n-k+1)\Ric ab \psi_{\ol abI} +
	    (n-k)(n-k+1) r\psi_I.
\end{align*}
Let $\phi\in\Gamma(\Wedge^{0,k}X\otimes K_X^{-k}\otimes V\otimes {V'}^{-1})$
be a solution of the equation
$(\db+\db^*)\phi=0$.
Then, since $\nabla\phi$ and $\nabla\psi$ are perpendicular,
 the contraction $\xi = \phi\psi\in\Gamma(K_X^{-k}\otimes V)$
 satisfies
$$
	\langle \nabla^2\xi, e_a e_{\ol a} \rangle
	 = (n-k)(n-k+1)r\xi - (n-k)\El\xi.
 $$
If we cancel $\El$ by using
$$
	\langle R\xi, e_a e_{\ol a} \rangle =
	  (n+1)(n-2k) r\xi - n \El\xi,
 $$
we obtain the following field equation.
\begin{defn}\label{defn:D0}
Let $\mathcal{D}_0$ be the differential operator on
$\Gamma(K_X^{-k}\otimes V)$:
$$
	\mathcal{D}_0 \xi =
	   \langle \nabla^2\xi, ke_ae_{\ol a} + (n-k)e_{\ol a}e_a \rangle
         - (n+2)k(n-k)r\xi.
 $$
\end{defn}
This can also be considered as a global operator on a conformally
Bochner-K\"ahler manifold. Then \eqref{eq:lap} is valid even
in the case $h=0$ if $\phi$ satisfies the field equation
$\mathcal{D}_0\phi=0$.

%
\section{The Penrose transform}\label{s:pt}

In this section we give a
generalization of the Penrose transform, which gives a relation
between the field equations given in the previous section and
cohomology groups of positive degree on twistor spaces.
\begin{lem}
We have an isomorphism
$$
	H^{k(n-k)}(G_{k,n}, {\cO}(H^{-n-h}))
	  \simeq \hS^h\Wedge^{0,k}\otimes (\Wedge^{n,0})^{-k}
 $$
as a representation space of $\Un(n)$.
\end{lem}
\begin{proof}
 If we consider the both spaces as $\SU(n)$-modules,
the equivariance is immediate from the Bott-Borel-Weil-Kostant
theorem (the BBWK-theorem).
Under this identification,
we can show easily that the  actions of a scalar matrix
coincide.
\end{proof}

As in the previous section, $X$ is an $n$-dimensional
conformally Bochner-K\"ahler manifold with $n>2$
and $V$ is a vector bundle on $X$ with
connection such that
$W=H^{-n-h}\otimes p^*V$ has a holomorphic structure.
We define the Penrose transform:
\begin{align*}
\begin{CD}
	  \PT : H^{k(n-k)}(Z_k(X), {\cO}(W))
	    @>>>
	  \Gamma(X,\Cup_{x\in X}
	    H^{k(n-k)}(Z_k(X)_x,{\cO}(W)))\\
	  @. @| \\
	  {} @. \Gamma(X, \hS^h\Wedge^{0,k}X\otimes K_X^{-k}\otimes V)
\end{CD}
\end{align*}
\begin{thm}\label{z2x}
If $h=0$, then an image of $\PT$ is a solution of the equation
$\mathcal{D}_0\phi = 0$.
If $h>0$,
an image of $\PT$ is
a solution of the equation $(\db + \db^*)\phi = 0$, that is,
 a harmonic form of the complex on~$X$:
\begin{align*}
\begin{split}
	\Gamma(\Wedge^{0,k-1}X\otimes \hS^{h-1} \Wedge^{0,k}X\otimes
	 K_X^{-k}\otimes V)& \\
	  \buildrel \db \over
	 \longrightarrow
	\Gamma(\hS^h \Wedge^{0,k}X\otimes &K_X^{-k}\otimes V)
	  \buildrel \db \over \longrightarrow\\
	 &\Gamma(\Wedge^{0,k+1}X\otimes \hS^{h-1} \Wedge^{0,k}X\otimes
	  K_X^{-k}\otimes V).
\end{split}
\end{align*}
Furthermore $\PT$ is injective.
\end{thm}

\begin{proof}
By the definition of the almost complex structure,
for $x\in X$,  $Z_k(X)_x$ is
a complex submanifold of $Z_k(X)$.
The normal bundle is a homogeneous vector bundle
$N \simeq N_1\oplus N_2$:
\begin{align*}
	N_1 = \SU(n)\times_{\kappa_1} \C^k,\qquad
	N_2 = \SU(n)\times_{\kappa_2} \C^{n-k},
\end{align*}
where $\kappa_1$ and $\kappa_2$ are
 representations of $S(\Un(k)\times\Un(n-k))$:
\begin{align*}
	\kappa_1(A,B) = \bar{A},\qquad
	\kappa_2(A,B) = B.
\end{align*}
Let $\kappa_3$ be the representation:
$$
	\kappa_3(A,B) = \det(A)^{-1}.
 $$
Then the associated line bundle is the hyperplane bundle.

Let $T=S(\Un(1)\times\dots\times\Un(1))$ be the subgroup of
diagonal matrices of $\SU(n)$.
By restricting the action of $\SU(n)$ or $S(\Un(k)\times\Un(n-k))$
to $T$, we consider weights of representations.
For an integer $i$ such that $1 \le i \le n-1$,
let $\lambda_i$ be the highest weight with respect to the
representation $\Wedge^i \C^n$.
As stated in \cite{[B.E]}, when we apply the BBWK-theorem,
it is convenient to consider the lowest weights of
irreducible representations of $S(\Un(k)\times\Un(n-k))$.
We can show, by simple computation, that the irreducible representations
$\Wedge^a \kappa_1$, $1\le a \le k$,
$\Wedge^b \kappa_2$, $1 \le b \le n-k$ and  $\kappa_3$
are characterized by their lowest weights
$-\lambda_a$, $-\lambda_{n-b}$ and $-\lambda_k$, respectively.

In order to prove that the Penrose transform is injective, let
$[\omega]$ be an element of $H^{k(n-k)}(\cO(W))$
such that the restriction $[\omega_x]\in
H^{k(n-k)}(Z_k(X)_x, \cO(W))$ vanishes for each $x\in X$.
Let $\omega_0 = \omega\in \Omega^{0,k(n-k)}(W)$ be
a representative of $[\omega]$.
Inductively,
for a non-negative integer $l$, let $\omega_l$ be
a representative of $[\omega]$
having the horizontal degree at least $l$.
We fix $x\in X$ for a while.
We have an isomorphism induced from the standard Hermitian metric
as a homogeneous vector bundle:
\begin{align*}
	\Wedge_H^{0,1} Z_k(X) \simeq {\bar N}_1^*\oplus {\bar N}_2^*
	 \simeq N_1 \oplus N_2,
\end{align*}
where $\Wedge_H^{0,1} Z_k(X)$ is the horizontal part of
$\Wedge^{0,1}Z_k(X)$.
For non-negative integers $a$ and $b$ such that $a+b=l$,
we can define an element
$$
	\omega_{a,b}(x) \in \Omega^{0,k(n-k)-l}_{G_{k,n}}(
		W\otimes \Wedge^a N_1\otimes \Wedge^b N_2)
 $$
by the canonical projection of $\omega_l$.
The $\db$-closedness of $\omega_l$ implies
that $\omega_{a,b}(x)$ is $\db$-closed.
If $l=0$, then by the assumption, there is a section
$$
        s_{0,0}(x)\in \Omega^{0, k(n-k)-1}_{G_{k,n}}(W)
 $$
such that $\db s_{0,0}(x) = \omega_{0,0}(x)$. By the BBWK-theorem,
we have
\begin{align*}
	H^{k(n-k)-1}(G_{k,n}, H^{-n-h}) &= 0.
\end{align*}
Hence $s_{0,0}(x)$ is unique and smooth with respect to $x$ by
the elliptic regularity.
In the case $l>0$,
by the BBWK-theorem, we have
\begin{align*}
	H^{k(n-k)-l-i}(G_{k,n}, H^{-n-h}\otimes \Wedge^a N_1
		\otimes \Wedge^b N_2) &= 0,\qquad i = 0,1.
\end{align*}
This means, for each $x\in X$, there exists a unique section
$$
	s_{a,b}(x)\in \Omega^{0, k(n-k)-l-1}_{G_{k,n}}(
		W\otimes \Wedge^a N_1\otimes \Wedge^b N_2)
 $$
such that $\db s_{a,b}(x) = \omega_{a,b}(x)$.
Again, by the elliptic regularity, $s_{a,b}(x)$ is smooth with respect to $x$.

Therefore if we consider $s_l = \sum_{a+b=l} s_{a,b}(x)$
as a section of $\Omega^{0,k(n-k)-1}(W)$,
 $\omega_l - \db s_l$ has horizontal degree at least $l+1$.
Hence, by induction on $l$,
 we have shown that the Penrose transform is injective.

Now we show that $\PT\alpha$ satisfies the field equation.
Put $Z=Z_k(X)$.
As stated in \cite{[I4]}, by a theorem of Kodaira (\cite{[Ko]}),
there exists a double fibration:
\begin{align*}
\setlength{\unitlength}{0.5mm}%
  \begin{picture}(50,30)
    \put(20,20){$Y$}
    \put(0,0){$Z$}
    \put(40,0){$X_\C$}
    \put(18,18){\vector(-1,-1){10}}
    \put(27,18){\vector(1,-1){10}}
    \put(7,16){$\scriptstyle p_1$}
    \put(34,16){$\scriptstyle p_2$}
  \end{picture}
\end{align*}
where $X_\C$ parametrizes submanifolds of $Z$
isomorphic to $G_{k,n}$ whose normal bundles are isomorphic to $N$.
The space $X_\C$ contains $X$ as the set of real fibers.
The manifold $Y$ is defined as:
$$
	Y = \{(z,x)\in Z\times X_\C \mid
 	  \text{$z$ is a point of the submanifold corresponding to $x$.}\}.
 $$
Hence $p_1\circ p_2^{-1}(x)$ is the submanifold corresponding
to $x\in X_\C$, and
the projection $p_2: Y \rightarrow X_\C$ is a holomorphic
fiber bundle with fiber $G_{k,n}$.
Moreover the tangent space $T_x X_\C$ is canonically isomorphic to
the space of global sections of the normal bundle of $p_1\circ p_2^{-1}(x)$.
For $x\in X$, let $Z_x = p_2^{-1}(x)$.
We compare two embeddings $Z_x\subset Y$ and $p_1(Z_x) \subset Z$.
Let $\cO_Y^l$ and $\cO_Z^l$ denote the $l$th order neighborhood
sheaves.
Then the $l$-jet of $\PT\alpha$ lies in the image of
\begin{align*}
\begin{CD}
	H^{k(n-k)}(p_1(Z_x), \cO_Z^l (W)) @>>>
	  H^{k(n-k)}(Z_x, \cO_Y^l (p_1^* W))\\
	@. @|\\
	{} @.  J_l(\hS^h\Wedge^{0,k}X\otimes K_X^{-k}\otimes V )_x
\end{CD}
\end{align*}
In general, let $B$ be a submanifold of $A$ with normal bundle $N$
and $E$ be a holomorphic vector
bundle over $A$. Then we have an exact sequence
of sheaves:
$$
	0 \longrightarrow
	  \cO_B(S^lN^*\otimes E) \longrightarrow
	  \cO_A^l(E) \longrightarrow
	  \cO_A^{l-1}(E) \longrightarrow 0.
 $$
Applying this sequence to $Z_x\subset Y$ and $p_1(Z_x)\subset Z$,
we have exact sequences:
\begin{multline*}
\begin{CD}
	0 @>>>
	H^{k(n-k)}(\cO(S^lN^*\otimes W)) @>>>
	  H^{k(n-k)}(\cO_Z^l(W))\\
	@. @VV{\sigma^*}V @VV{p_1^*}V\\
	0 @>>>
	(S^lT^*X\otimes \hS^h \Wedge^{0,k}X\otimes K_X^{-k}\otimes V)_x @>>>
	  J_l(\hS^h\Wedge^{0,k}X\otimes K_X^{-k}\otimes V)_x
\end{CD}\\[4mm]
\begin{CD}
	 @>>> H^{k(n-k)}(\cO_Z^{l-1}(W)) @>>> 0\\
	 @.    @VV{p_1^*}V \\
	 @>>> J_{l-1}(\hS^h\Wedge^{0,k}X\otimes K_X^{-k}\otimes V)_x @>>> 0
\end{CD}
\end{multline*}
where $H^{k(n-k)-1}$ terms disappear because of the BBWK-theorem.

  If $l=0$, $p_1^*$ is an isomorphism by definition.
For all $l$ the map $\sigma^*$ is injective as can be shown by
the BBWK-theorem. Hence by induction $p_1^*$  is injective
for all $l$ and maps
$H^{k(n-k)}(\cO_Z^l(W))$ into some subspace of
$J_l(\hS^h\Wedge^{0,k}X\otimes K_X^{-k}\otimes V)_x$.

  If $h \ge 1$ and $l = 1$, then the $1$-jet bundle does not depend
on the Hermitian metric of $X$. Hence, by calculating in the
flat case,
we see that this
is the  subspace corresponding to the
equation $(\db + \db^*)\PT\alpha = 0$.

  When $h=0$, $p_1^*$ gives an isomorphism so $\PT\alpha$ satisfies no
first order equation. Passing to the second order
neighborhood, however, we do obtain a proper subspace of
$J_2(K_X^{-k}\otimes V)$.
Thus $\PT\alpha$ satisfies a second order equation $\mathcal{D}_*\PT\alpha=0$.
By definition, we can assume that
the highest order part of $\mathcal{D}_*$ is the same as
for $\mathcal{D}_0$.
Since the problem is local, we can assume that there is a
line bundle $V'$ on $X$ such that $H\otimes p^*V'$ is
a holomorphic line bundle.
We consider the natural product map
\begin{multline*}
	H^0(p_1(Z_x),\cO^2_Z(H\otimes p^*V'))\otimes
	  H^{k(n-k)}(p_1(Z_x),\cO^2_Z(W\otimes H^{-1}\otimes p^*{V'}^{-1}))\\
	\longrightarrow
	  H^{k(n-k)}(p_1(Z_x),\cO^2_Z(W))
	\longrightarrow
	  J_2(K_X^{-k}\otimes V)_x.
\end{multline*}
Here every holomorphic section of
$H\otimes p^*V'$ on the first order neighborhood of $p_1(Z_x)\subset Z$
extends uniquely to the second order neighborhood by the
BBWK-theorem. Hence the image is non-zero and, by the
computation in the previous section, it is
annihilated by both $\mathcal{D}_*$ and $\mathcal{D}_0$.
  Thus $\mathcal{D}_*-\mathcal{D}_0$ is a first order
operator whose homomorphism annihilates the image of
$$
	H^0(p_1(Z_x),\cO^1_Z(H\otimes p^*V'))\otimes
          H^{k(n-k)}(p_1(Z_x),\cO^1_Z(W\otimes H^{-1}\otimes p^*{V'}^{-1}))
 $$
in $J_1(K_X^{-k}\otimes V)_x$.
Again, since the $1$-jet bundle is not intervened by the Hermitian
metric of $X$, by computing in the flat case,
this is shown to be the whole space.
Hence we conclude that $\mathcal{D}_* = \mathcal{D}_0$.

Therefore, for every $h\ge 0$, $\PT\alpha$ satisfies the desired equation.
\end{proof}

%
\section{The inverse Penrose transform}\label{sec:inv}

In this section, we prove that the Penrose transform in the previous
section is surjective.
This is done by
constructing explicitly a Dolbeault representative
corresponding to a solution of the field equation on the base manifold.

The idea is essentially same as in the case of Riemannian manifolds.
But the presence of the curvature tensor makes the
proof of well-definedness
far more complicated,
and many intermediate operators appear during the computation.
For convenience, we summarize them in Appendix.

Since the problem is local, we can assume that $X$ is
an $n$-dimensional Bochner-K\"ahler manifold with $n>2$.
Let $V$ be a vector bundle on $X$ with connection
satisfying the condition
that $H^{-n-h}\otimes p^*V$ is a holomorphic vector bundle.

First, we give a lemma about
 the curvature tensor of a Bochner-K\"ahler manifold.
Let $(e_a)_a$ be a local orthonormal frame of the holomorphic tangent
bundle of $X$. We use notation of
curvature tensors defined in \S\ref{sec:twh}.

By a theorem of Kamishima (\cite{[Ka]}), a Bochner-K\"ahler manifold
is locally isomorphic to (i) $\C^n$ with a flat metric,
  or (ii) $\Hypb^{l}\times  \Proj^{n-l}$
for some integer $l$ such that $0 \le l \le n$.
This means that the Bochner-flat condition also strongly
 restricts the Ricci tensor.
In this paper, we only need the following lemma, which is
deduced immediately from that theorem.

\begin{lem}\label{curvtoric}
The Ricci tensor of a Bochner-K\"ahler manifold
 $X$ is covariantly constant and
 satisfies
$$
	\Ric ab \Ric bc = r \Ric ac + s \delta_a^c,
 $$
where $s$ is a scalar function.
\end{lem}

Let $\nabla e_a = \omega_a^be_b$ be the
connection form of the Levi-Civita connection.
Since a connection defines horizontal lifts of vector fields on $X$
to those on $Z_k(X)$, $e_a$ is also considered to be a vector field
on $Z_k(X)$.
By restricting the $\SO(2n)$-action on $Z$, we have
a $\Un(n)$-action on $G_{k,n}$.
It is convenient however to consider the
complexified $\GL(n;\C)$-action on $G_{k,n}$. It defines a
$\C$-linear map:
\begin{align*}
	\FF: \gl(n;\C) \longrightarrow H^0(G_{k,n},
		\cO(T^{(1,0)}G_{k,n})).
\end{align*}
With respect to the Lie algebra structures, we have:
\begin{equation}\label{eq:lie_alg}
	\FF([a,b]) = - [\FF(a),\FF(b)].
\end{equation}
Let $(E^a_b)_{a,b}$ denote the standard basis of $\gl(n;\C)$ such that
$E^a_b e_c = \delta^a_c e_b$.
We define
\begin{equation}\label{eq:vf}
	\FF^a_b = - \FF(\E^a_b).
\end{equation}

Before giving an explicit description of $\FF^a_b$,
we give here a few relations between the Pl\"uker coordinates
induced from the Pl\"uker relation.

\begin{lem}\label{lem:defeq}
  \begin{enumerate}
    \item Let $J$, $K$ be multi-indices of length $k$.
	Then we have, on $G_{k,n}$
\begin{align*}
	Z^{a\ol cJ}Z^{\ol bcK} =
 	  - Z^JZ^{\ol baK} + Z^{\ol baJ} Z^K.
\end{align*}
    \item For a multi-index $J$ of length $k$, we have, on $G_{k,n}$
\begin{align*}
	Z^{\ol ab\ol cdJ}Z^J =
	   Z^{\ol abJ}Z^{\ol cdJ} + Z^{\ol adJ}Z^{b\ol cJ}.
\end{align*}
\end{enumerate}
\end{lem}
\begin{proof}
  When $a\ne b$, by using \eqref{eq:plucker}, we compute:
\begin{align*}
	Z^{a\ol cJ}Z^{\ol bcK}
	  &= Z^{a\ol cJ}Z^{\ol bcK} - Z^{\ol caJ}Z^{c\ol bK}\\
	  &= Z^{a\ol a J} Z^{\ol baK} - Z^{\ol a aJ}Z^{a\ol bK}
	     + Z^{a\ol bJ}Z^{\ol bbK} - Z^{\ol baJ}Z^{b\ol bK}\\
	  &= -Z^JZ^{\ol baK} + Z^{\ol baJ} Z^K.
\end{align*}
The case $a=b$ can be proved in the same way and we have proved (1).

By putting $J=\ol dc K$ in (1), we have (2).
\end{proof}

For a multi-index $I$ of length $k$, we take local fiber coordinates
as in \S\ref{sec:twh}. That is,
$$
	w_{ij}=Z^{\ol \imath jI}/Z^I,\qquad i\in I,\quad j\not\in I
 $$
are local coordinates on
$$
	U_I=\{(Z^J)_{J}\in G_{k,n} \mid Z^I\ne 0\}
 $$
with the center $z_0$. For a general multi-index $J$ of length $k$,
put $z^J = Z^J/Z^I$.

\begin{lem}\label{vectorfield}
The vector field $\FF^a_b$ is written in the local coordinates as:
\begin{align*}
	\FF^a_b = - z^{a\ol \imath I}z^{\ol bjI}\partdel{}{w_{ij}}.
\end{align*}
\end{lem}

\begin{proof}
	For $a\not=b$, we compute the infinitesimal
action of the one-parameter subgroup:
$$
	A(t) = \exp(tE^a_b).
 $$
Let $(e_J)_J$ be the standard basis of $\Wedge^k \Cn$.
Then the action of the matrix is written as:
$$
	A(t)e_J = e_J + t e_{\ol abJ}.
 $$
Hence the action on a Pl\"ucker coordinate is written as
$$
	Z^J(t) = Z^J + tZ^{\ol baJ}.
 $$
Therefore, by using Lemma \ref{lem:defeq} (2), we have
\begin{align*}
	\left. \frac{d}{dt}w_{ij}(t) \right|_{t=0}
	  = z^{\ol ba\ol \imath jI} - z^{\ol baI} z^{\ol \imath jI}
	= z^{a\ol \imath I} z^{\ol bjI},
\end{align*}
 for $i\in I$ and $j\not\in I$.

  When $a=b$, let
$$
	B(t) = \exp(t\im E^a_a),
 $$
where we remark that the index $a$ is fixed
 and Einstein's convention for summation is not applied.
Then, in a similar way, we compute
\begin{align*}
	\left. \frac{d}{dt}w_{ij}(t) \right|_{t=0}
	   = \im z^{a\ol \imath I} z^{\ol a jI}.
\end{align*}
Hence we complete the proof.
\end{proof}

Let $(e^a)_a$ be the dual frame of $\Wedge^{1,0}X$.
As in the case of $e_{\ol a}$, let $e^{\ol a}$,
$\FF^{\ol a}_{\ol b}$,  $R^{\ol a}_{\ol b}$, and $\omega^{\ol a}_{\ol b}$ be
the complex conjugate of
$e^a$, $\FF^a_b$, $R^a_b$ and $\omega^a_b$, respectively.
Put $W = H^{-n-h}\otimes p^*V$.
Now we define operators acting on
$\Omega^*(Z_k(X), W)$:
\begin{align}
	T_0 &= - i(\FF^b_a)R^a_b - i(\FF^{\ol b}_{\ol a})R^{\ol a}_{\ol b},
		\label{eq:T0}\\
	\La a &= L_{e_a} - \omega^b_a i(e_b)
		- [i(e_a),T_0],\label{eq:La}\\
	\Lb a &= L_{e_{\ol a}} -
		 \omega^{\ol b}_{\ol a} i(e_{\ol b})
		- [i(e_{\ol a}),T_0],\label{eq:Lb}\\
	D^{a} &= e^b i(\FF^{\ol b}_{\ol a}),\label{eq:D^a}\\
	D^{\ol{a}} &= - e^{\ol b} i(\FF^{\ol a}_{\ol b}),
		\label{eq:D^b}
\end{align}
where $L$ denotes the Lie derivative. Vector fields and forms on $X$
are considered to be those on $Z_k(X)$ in a natural way.
Note that $D^{\ol a}$ is not the complex conjugate of $D^a$.

\begin{lem}\label{lem:gr}
\begin{enumerate}
\item
Let $e_a' = e_b h^b_a$ be another local orthonormal frame of $T^{(1,0)}X$.
Let ${\La a}'$, ${\Lb a}'$,  ${D'}^{a}$ and ${D'}^{\ol{a}}$
be operators defined as above with respect to the frame $(e_a')_a$.
Then we have
$$
	{\La a}' = h^b_a \La{b},\quad
	{\Lb a}' = \overline{h^b_a}\Lb b,\quad
	 {D'}^a = D^b(h^{-1})^a_b ,\quad
	 {D'}^{\ol{a}} = D^{\ol{b}}\overline{(h^{-1})^a_b}.
 $$
\item Operators $\La a$, $\Lb a$, $D^{a}$ and $D^{\ol a}$ preserve
the double grading as differential forms.
That is, for non-negative integers $l$ and $l'$,
they map $\Omega^{l,l'}(W)$ to itself.
\end{enumerate}
\end{lem}
\begin{proof}
The transformation rules in (1) are immediate by
\begin{align*}
	L_{fv} &= f L_v  + df \wedge i(v),\\
	{\omega'}^a_b &= (h^{-1})^a_c{dh}^c_b + (h^{-1})^a_c\omega^c_d h^d_b,
\end{align*}
where $v$ is a vector field and $f$ is a function.

The operators $D^a$ and $D^{\ol a}$ preserve the grading
because we can write them as:
\begin{align*}
	D^a &= - z^{a\ol d\ol{I}}
		 (z^{\ol b c\ol{I}} e^b)
		 i(\partdel{}{{\ol w}_{cd}}),\\
	D^{\ol a} &= z^{\ol a c\ol{I}}
		 (z^{b\ol d\ol{I}} e^{\ol b})
		 i(\partdel{}{{\ol w}_{cd}}).
 \end{align*}

We can assume that the connection forms of $T^{(1,0)}X$
and $V$ vanish at $x_0$ and
we compute values at $(x_0,z_0)$.
By \eqref{eq:ddwa}, we have
\begin{align*}
	[L_{e_a}, \dwa] &= [i(e_a), R^j_i],\\
	[L_{e_a}, \dwb] &= [i(e_a), R^{\ol \jmath}_{\ol \imath}],
 \end{align*}
where $\dwb$ is the complex conjugate of $\dwa$.
On the other hand, we have
\begin{align*}
	[- \omega^b_a i(e_b) - [i(e_a),T_0], \dwa] &= - [i(e_a), R^j_i],\\
	[- \omega^b_a i(e_b) - [i(e_a),T_0], \dwb] &=
		- [i(e_a), R^{\ol \jmath}_{\ol \imath}].
 \end{align*}
Consequently we have
\begin{align*}
	[\La a, \dwa] = [ \La a, \dwb] = 0.
 \end{align*}
With respect to horizontal forms, we compute
\begin{align*}
	[\La a, z^{\ol b \ol{J}} e^b ] &=
	[\La a, z^{bJ} e^{\ol b} ] =
		0,\qquad \abs{J}=k-1,\\
	[\La a, z^{b \ol{J}} e^{\ol b}] &=
	[\La a, z^{\ol b J} e^b] =
		 0,\qquad \abs{J}=k+1.
 \end{align*}
For a function $f$, we have
$$
	[\La a, f] = e_a(f).
 $$
Hence $\La a$ preserves the grading. The proof for $\Lb a$ can
be done in the same way.
\end{proof}

Now we define operators which are used to construct
the inverse Penrose transform.
\begin{align}
	\Da &=  D^{a} \La a,\label{eq:Da}\\
	\Db &= D^{\ol{a}}\Lb a ,\label{eq:Db}\\
	\Dr &= - \Ric{a}{b}D^{\ol{b}}D^a.\label{eq:Dr}
\end{align}
These operators are globally well-defined
and preserve the grading of forms by the above lemma.

Put
$$
	F(x) = \sum_{i=0}^{\infty} \frac{x^i}{(i!)^2}.
 $$
This function and its derivatives play an important role by the
following property.
\begin{lem}\label{lem:diff_eq}
Let $l$ be a non-negative integer. Then the $l$th derivative of $F$
satisfies
\begin{align*}
	xF^{(l+2)}(x) + (l+1)F^{(l+1)}(x) - F^{(l)}(x) = 0.
\end{align*}
\end{lem}

	Let $\Wedge_V^{0,k(n-k)}$ denote the line subbundle of
$\Wedge^{0,k(n-k)} Z_k(X)$ spanned by vertical forms,
which are defined by the Levi-Civita connection.
If we identify $H^{-1}$ with $\bar H$ by the standard Hermitian metric,
we have
$$
	\Wedge^{0,k(n-k)}Z_k(X)\otimes H^{-n-h} \supset
		\Wedge_V^{0,k(n-k)}\otimes H^{-n-h}
		\simeq (H^{h}\otimes K_X^{k})^*,
 $$
where the canonical bundle appears because we consider
the action of $\Un(n)$ rather than $\SU(n)$.
On the other hand, we have an isomorphism
$H^0(G_{k,n},{\cO}(H^h)) \simeq \hS^h\Wedge^{k,0}X$
by the BBWK-theorem.
Hence we obtain
$$
	j: \Gamma(X, \hS^h\Wedge^{0,k}X\otimes
		K_X^{-k}\otimes V) \longrightarrow
		\Omega^{0,k(n-k)}
		  (Z_k(X), W).
 $$
We put
\begin{align*}
	f(x) &= (k+h-1)! F^{(k+h-1)}(x),\\
	g(x) &= (n-k+h-1)! F^{(n-k+h-1)}(x).
\end{align*}

\begin{defn}\label{defn:inv}
We define
$$
\begin{matrix}
	\AAA:& \Gamma(X, \hS^h\Wedge^{0,k}X\otimes
		K_X^{-k}\otimes V) &\longrightarrow
		&\Omega^{0,k(n-k)}(Z_k(X), W)\\[3mm]
		& \phi &\longmapsto
		& f(D_\beta + \Dr) g(D_\alpha)
			\jjj.
\end{matrix}
 $$
\end{defn}

\begin{rem}
This is well-defined since the differential operators
in the formula are even operators and $\Db$ and $\Dr$ are mutually
commutative.
By the proof of
 Lemma \ref{lem:gr}, the operators $\Da^i$, $\Db^i$ and $\Dr^i$
vanish for sufficiently large $i$.
\end{rem}

\begin{thm}\label{thm_x2z}
The form $\AAA(\phi)$ is $\db$-closed if $\phi$ satisfies
the differential equation stated in Theorem \ref{z2x}.
The restriction of $\AAA$ to the space of harmonic forms gives
 the inverse of the Penrose transform.
\end{thm}

The remainder of this section is devoted to a proof of this theorem.
Since the operators $\Da$, $\Db$ and $\Dr$ decrease the vertical grading
of forms, the last statement of the theorem follows immediately
if we show the first statement.

 Let $(e_a)_a$ be a local frame of $T^{(1,0)}X$ such  that
the connection form of the Levi-Civita connection vanishes at $x_0$.
Let $(e^{\ol I})_{\ol I}$ denote the associated frame of $\Wedge^{0,k}X$.
Instead of determining an explicit basis of $\hS^h\Wedge^{0,k}X$,
we treat it as a subbundle of $\Otimes^h\Wedge^{0,k}X$.
On the other hand, we consider an element of $\Otimes^h\Wedge^{0,k}X$
as an element of $\hS^h\Wedge^{0,k}X$ by the canonical projection.
Let $I_1,\dots,I_h$ be multi-indices
of length $k$.
Put
$$
	e^{\ol{I}_1,\dots,\ol{I}_h} =
	  e^{\ol{I}_1}\otimes\dots\otimes e^{\ol{I}_h},
 $$
which is considered as an element of $\hS^h\Wedge^{0,k}X$ as mentioned above.
Since we have taken a local frame of $T^{(1,0)}X$, we have a canonical
trivialization of $K_X$.
We also take a local frame of $V$ such that the connection form
vanishes at $x_0$.
Since the choice of the local frame of $V$
 does not affect the appearance of the computation,
we omit its index.
In fact, the information of $V$ which is needed for computation
is already encoded in \eqref{eq:spin_curv} and \eqref{eq:lap}.
Thus we can consider $e^\iii$ as a section of
$\hS^h\Wedge^{0,k}X\otimes K_X^{-k}\otimes V$.
Then we put
$$
	s^\iii = j(e^{{\ol I_1},\dots,{\ol I_h}}).
 $$
Define
\begin{align}
	E^{a} &= [\db, D^{a}] + \omega^a_bD^b,
	  \label{eq:E^b}\\
	E^{\ol a} &= [\db, D^{\ol a}] + \omega^{\ol a}_{\ol b}D^{\ol b}.
	  \label{eq:E^a}\
\end{align}
They satisfy the same transformation rules as $D^a$ and $D^{\ol a}$,
 respectively.

\begin{lem}\label{dirty}
For local coordinates $(w_{ij})_{i,j}$,
we have
      \begin{align*}
	\frac{\partial z^J}{\partial w_{ij}} = -z^{i\ol \jmath J}.
      \end{align*}
\end{lem}
\begin{proof}
  It is obvious
 if $|J\setminus I| \le  1$. If $|J\setminus I| \ge 2$,
    let $i_k\in I\setminus J$, $j_k\in J\setminus I$, $k = 1,2$, be numbers
   such that
    $i_1 \ne i_2$, $j_1 \ne j_2$. Then, by Lemma \ref{lem:defeq} (2),
 we have
    \begin{align*}
      Z^J = \frac{1}{Z^{\ol{\jmath}_1i_1\ol{\jmath}_2i_2J}}
	(Z^{\ol{\jmath}_1i_1J}Z^{\ol{\jmath}_2i_2J}
	 + Z^{{\ol\jmath}_1i_2J}Z^{i_1\ol{\jmath}_2J}).
    \end{align*}
Hence we obtain inductively the desired formula.
\end{proof}

\begin{lem}\label{lem:act_e}
We have
\begin{align*}
	E^{a} s^\iii &\equiv
	     - (n-k+h)e^a s^\iii
	    - \sum_{j=1}^h
		e^b s^{\ol{I}_1,\dots,\ol{b}a\ol{I}_j,\dots,\ol{I}_h},\\
	E^{\ol a} s^\iii &\equiv
	    -(k+h)e^{\ol a} s^\iii
	    -\sum_{j=1}^h
		e^{\ol b}
		s^{\ol{I}_1,\dots,b\ol a\ol{I}_j,\dots,\ol{I}_h},
\end{align*}
where we consider equivalence modulo $(1,0)$-forms.
\end{lem}
\begin{proof}
Let $\rho^{\ol I}$ be the image of the standard trivialization
of $\bar{H}$
by the isomorphism $\bar{H}\simeq H^{-1}$.
Let $K^{\ol I}$ be the standard trivialization of $\Wedge^{0,k(n-k)}_V$.
 Then, by definition, we have
\begin{equation}\label{eq:siii}
	s^\iii =
	 \rho^{\ol{I}_1}\otimes\dots\otimes \rho^{\ol{I}_h}\otimes
	 {\rho^{\ol I}}^{\otimes(n)}\otimes K^{\ol I}.
\end{equation}
Put
$$
	N = \sum_J \left|\frac{z^J}{z^I}\right|^2.
 $$
First, we compute:
\begin{align*}
	L_{\FF^{\ol a}_{\ol b}} \rho^{\ol I} &=
		\nabla_{\FF^{\ol a}_{\ol b}} \rho^{\ol I}
		= -\frac{\FF^{\ol a}_{\ol b}(N)}{N} \rho^{\ol I}\\
		&= - \frac{z^J z^{\ol\imath j\ol{J}} z^{\ol ai\ol I}
				z^{b\ol\jmath\ol{I}}}
			     {N}
			\rho^{\ol I}
			&&\text{[By Lemma \ref{dirty}]}\\
		&= - 
			\frac{z^J z^{\ol ai\ol I}z^{b\ol\imath \ol{J}}}{N}
			\rho^{\ol I}
			&&\text{[By Lemma \ref{lem:defeq} (1)]}\\
		&= ( 
			\frac{z^J z^{b\ol a\ol J}}{N}
			- z^{b\ol a\ol I})
			\rho^{\ol I}
			&&\text{[By Lemma \ref{lem:defeq} (1)]}
\end{align*}
By changing indices, we obtain:
\begin{align*}
	L_{\FF^{\ol a}_{\ol b}} \rho^{\ol I} =
		(\delta^b_a + z^{\ol{a}b\ol I}
		 + \frac{z^{\ol{a}bJ} z^{\ol J}}{N})
		\rho^{\ol I}.
\end{align*}
Second, we compute:
\begin{align*}
	L_{\FF^{\ol a}_{\ol b}} d{\ol w}_{ij} =
		d\FF^{\ol a}_{\ol b}(\ol{w}_{ij}) =
		(-\delta_i^b z^{\ol ai\ol I} +
		  \delta_a^j z^{b\ol\jmath\ol{I}})d\ol{w}_{ij} + \dots
\end{align*}
Hence we have
\begin{align*}
	L_{\FF^{\ol a}_{\ol b}} K^{\ol I} 
	= (-(n-k) z^{\ol ab\ol I} +
	     k z^{b\ol a \ol I} ) K^{\ol I}
	= - (k\delta_a^b + n z^{\ol{a}b\ol I})
			K^{\ol I}.
\end{align*}
Thus, we have
\begin{align*}
	L_{\FF^{\ol a}_{\ol b}} s^\iii =
	     (n-k+h)\delta_a^b s^\iii
	    + \sum_{j=1}^h
		  s^{\ol{I}_1,\dots,\ol ab\ol{I}_j,\dots,\ol{I}_h}
	    +(n+h)
		\frac{z^{\ol abJ} z^{\ol J}}{N} s^\iii.
\end{align*}
Consequently we have the first equality.
By changing indices we have:
\begin{align*}
	L_{\FF^{\ol a}_{\ol b}} s^\iii =
	    -(k+h)\delta_a^b s^\iii
	    - \sum_{j=1}^h
		  s^{\ol{I}_1,\dots,b\ol a\ol{I}_j,\dots,\ol{I}_h}
	    -(n+h)
		\frac{z^{b\ol{a}J} z^{\ol J}}{N} s^\iii,
\end{align*}
from which the second equality follows.
\end{proof}

\begin{lem}\label{lem:com_ed}
We have
\begin{align*}
	[E^a, D^b] &= - e^a D^b - e^b D^a,\\
	[E^{\ol a}, D^{\ol b}] &= - e^{\ol a}D^{\ol b}
			- e^{\ol b}D^{\ol a},\\
	[E^{\ol a}, D^b] &= \delta_a^b B_0,
 \end{align*}
where
\begin{equation}\label{eq:B0}
	B_0 = e^a D^{\ol a} = D^{a} e^{\ol a}.
\end{equation}
\end{lem}

\begin{proof}
These formulas follow from
$$
	[\FF^a_b, \FF^c_d] = \delta^a_d \FF^c_b - \delta^c_b \FF^a_d,
 $$
which is an immediate consequence of \eqref{eq:lie_alg}.
\end{proof}

We define
\begin{align}
	\dda &= e^a \La a,
		\label{eq:da}\\
	\ddb &= e^{\ol a} \Lb a,
		\label{eq:db}\\
	\ddv &= - i(\FF^{\ol b}_{\ol a})R^{\ol a}_{\ol b}.
		\label{eq:dv}
\end{align}
Obviously these operators do not preserve the grading
of forms. But we see, by Lemma \ref{lem:gr},
 they do not decrease the first
grading. So we can consider them as operators acting on
$$
	\frac{\Omega^*(Z_k(X), W)}
	   {\Oplus_{i=1}^{\infty} \Omega^{i,*}(Z_k(X),W)} \simeq
	\Omega^{0,*}(Z_k(X), W).
 $$

Now we compute $\db\AAA(\phi)$.
We split it into two steps
\begin{align}
	\db \gj &= (\ddb + \ddv)\gj,
		\tag{\Romannumeral 1} \label{eq:dbxg}\\
	[\db, f(\Db+\Dr)] \gj
		&= - f(\Db+\Dr) (\ddb + \ddv) \gj,
		\tag{\Romannumeral 2} \label{eq:dbfg}
\end{align}
from which the theorem follows immediately.

\begin{proof}[Proof of \eqref{eq:dbxg}]
Let $\phi = \phi_\iii e^\iii$. Then we have
$$
	\jjj = \phi_\iii s^\iii
 $$
and
$$
	\db \jjj \equiv d \phi_\iii \wedge s^\iii
		+ \phi_\iii d s^\iii.
 $$
The value of the first term at $(x_0,z_0)$ is
$$
	d \phi_\iii \wedge s^\iii = ( \dda + \ddb )\jjj.
 $$
By \eqref{eq:ddwa}, we have
$$
	d K^{\ol I} |_{(x_0,z_0)} = \ddv K^{\ol I}.
 $$
Hence, by \eqref{eq:siii}, we have
\begin{align*}
	\db\jjj = ( \dda + \ddb + \ddv ) \jjj.
\end{align*}

By an inductive argument, we have
$$
	\Da^i\jjj = D^{a_1}\dots D^{a_i}
		 j\langle
		   \nabla^i\phi ,
		   e_{a_i}\dots e_{a_1}
		 \rangle.
 $$
Hence we compute
$$
	\db \Da^i\jjj = [\db,D^{a_1}\dots D^{a_i}]
                 j\langle
                   \nabla^i\phi ,
                   e_{a_i}\dots e_{a_1}
                 \rangle
	  + (\dda + \ddb + \ddv) \Da^i\jjj.
 $$
Therefore  we obtain
$$
	(\db - \ddb - \ddv)\Da^i\jjj =
		[\db,D^{a_1}\dots D^{a_i}]
                 j\langle
                   \nabla^i\phi ,
                   e_{a_i}\dots e_{a_1}
                 \rangle
	  + \dda \Da^i\jjj.
 $$
In order to compute the first term of the right-hand side, put
\begin{equation}\label{eq:Ea}
	\Ea = E^{a}\La a.
\end{equation}
Then, by Lemma \ref{lem:act_e} and the assumption on $\phi$, we have
$$
	\Ea\jjj = - (n-k+h)\dda\jjj.
 $$
By Lemma \ref{lem:com_ed}, we compute
\begin{align*}
	[\Ea,\Da]\Da^i\jjj = -2 \dda \Da^{i+1}\jjj,
\end{align*}
where we use
$$
	[\dda ,\Da]\Da^i\jjj = 0,
 $$
which follows from the fact that
the curvature form is of type $(1,1)$ with respect to
the original complex structure of $X$.
This also means
$$
	[[\Ea,\Da],\Da]\Da^i\jjj = 0.
 $$
Hence we obtain
$$
	\Ea\Da^i\jjj = -2i\dda\Da^i\jjj - (n-k+h) \dda \Da^i\jjj.
 $$
Therefore we compute
\begin{align*}
	[\db,D^{a_1}\dots D^{a_i}]
             j\langle
               \nabla^i\phi ,
                  e_{a_i}\dots e_{a_1}
             \rangle
	&= \sum_{j=1}^i \Da^{j-1} \Ea \Da^{i-j}\jjj\\
	&= -i(i-1)\dda\Da^{i-1}\jjj - i (n-k+h) \dda \Da^{i-1}\jjj.
\end{align*}
Thus, by Lemma \ref{lem:diff_eq}, we have
\begin{align*}
	(\db - \ddb - \ddv) \gj
	&= \dda (- g''(\Da)\Da - (n-k+h) g'(\Da) + g(\Da)) \jjj\\
	& = 0,
\end{align*}
which completes the proof.
\end{proof}

\begin{proof}[Proof of \eqref{eq:dbfg}]
We have
$$
	\Db^i\Dr^j \jjj =
	  \Dr^j D^{\ol{a}_1}\dots D^{\ol{a}_i}
	  j\langle
	    \nabla^i\phi, e_{\ol{a}_i}\dots e_{\ol{a}_1}
	  \rangle.
 $$
Hence, by putting
\begin{align}
	\Eb &= E^{\ol{a}} \Lb a,
		\label{eq:Eb}\\
	\Er &= [\db, \Dr],
		\label{eq:Er}
\end{align}
we obtain
\begin{align*}
	[\db, \Db+\Dr ] (\Db+\Dr)^i \gj
	 =
		(\Eb + [\dda,\Db]+\Er) (\Db+\Dr)^i \gj.
\end{align*}

Now, as in the proof of \eqref{eq:dbxg},
we compute the commutation relation between $\Eb + [\dda,\Db]+\Er$
and $\Db+\Dr$.
By Lemma \ref{lem:com_ed}, we have
$$
	[\Eb, \Db](\Db+\Dr)^i\gj = - 2 \Db\ddb (\Db+\Dr)^i\gj
 $$
and
\begin{align*}
	[\Eb, \Dr] (\Db+\Dr)^i\gj
	   = (\Db B_1 - \Dr \ddb - B_0 C_1) (\Db+\Dr)^{i} \gj,
\end{align*}
where
\begin{align}
	B_1 &=  \Ric ab D^a e^{\ol b},
		\label{eq:B1}\\
	C_1 &= \Ric ab D^{\ol b} \Lb a.
		\label{eq:C1}
\end{align}
We compute
\begin{multline*}
	[\dda,\Db]\Db^{i+1}\Da^j \jjj\\
\begin{split}
	&= e^aD^{\ol b} D^{\ol c}
	   D^{{\ol c}_1}\dots D^{{\ol c}_i}
	   D^{d_1}\dots D^{d_j}
	  j\langle
		R\nabla^{i+j+1}\phi, e_{d_j}\dots e_{d_1}
		  e_{\ol{c}_i}\dots e_{\ol{c}_1}
		  e_{\ol c}e_{\ol b} e_a
	  \rangle \\
	&= e^aD^{\ol b} D^{\ol c}
	   D^{{\ol c}_1}\dots D^{{\ol c}_i}
	   D^{d_1}\dots D^{d_j}
	  j\langle
		- R^{\ol d}_{\ol c}\nabla_{e_{\ol d}}\nabla^{i+j}\phi,
		  e_{d_j}\dots e_{d_1}
		  e_{\ol{c}_i}\dots e_{\ol{c}_1}
		  e_{\ol b} e_a
	  \rangle\\
	 &\quad + \Db[\dda,\Db]\Db^i\Da^j \jjj.
\end{split}
\end{multline*}
Obviously we have
\begin{align*}
	[[\dda,\Db],\Dr] \Db^i\Dr^j\gj = 0.
\end{align*}
Therefore we obtain
\begin{align*}
	[[\dda,\Db],\Db] (\Db+\Dr)^i\gj &= 2 T_1 (\Db+\Dr)^i\gj,\\
	[[\dda,\Db],\Dr] (\Db+\Dr)^i\gj &= 0,
\end{align*}
where
\begin{equation}\label{eq:T1}
	T_1 = \Curv abcd e^a D^{\ol b} D^{\ol d} \Lb c.
\end{equation}
We compute
\begin{multline*}
	[\Er,\Db](\Db+\Dr)^i\gj\\
  \begin{split}
	&=  -
	    \Ric ab (D^a [E^{\ol b},\Db] + D^{\ol b} [E^{a},\Db])
		(\Db+\Dr)^i\gj\\
	&= (\Db B_1 - \Dr \ddb -B_0 C_1)
		(\Db+\Dr)^i\gj,
  \end{split}
\end{multline*}
and
\begin{multline*}
	[\Er,\Dr](\Db+\Dr)^i \gj\\
  \begin{split}
	&= -
	  \Ric ab (D^a [E^{\ol b},\Dr] + D^{\ol b}[E^a,\Dr])(\Db+\Dr)^i\gj\\
	&= 2\Dr( B_1 + B_2 - rB_0)
	  (\Db+\Dr)^i \gj\\
	&= -2 \Dr\ddv (\Db+\Dr)^i \gj,
  \end{split}
\end{multline*}
where we put
\begin{equation}\label{eq:B2}
	B_2 =  \Ric ab e^a D^{\ol b},
\end{equation}
and we use Lemma \ref{curvtoric} and
\begin{align}
	D^{\ol a} D^a &= 0,\notag\\
	B_1 + B_2 - r B_0  &= -\ddv.
	  \label{eq:ddv2B}
\end{align}
The first equation is immediately obtained by \eqref{eq:plucker}.
Therefore, by using
$$
	T_1 + \Db B_1 - B_0 C_1 
	 =  -\Db \ddv,
 $$
we obtain
\begin{multline*}
	[[\db, \Db+\Dr], \Db+\Dr] (\Db+\Dr)^i \gj\\
	 =
		 - 2 (\ddb+\ddv) (\Db+\Dr)^{i+1} \gj.
\end{multline*}
Since
$$
	[[[\db, \Db+\Dr], \Db+\Dr], \Db+\Dr] (\Db+\Dr)^i \gj = 0,
 $$
by using Lemma \ref{lem:diff_eq}, we finally obtain
\begin{align*}
	[\db,f(\Db+\Dr)]\gj
		 &= f'(\Db+\Dr)(\Eb + [\dda, \Db] + \Er)\gj \\
		 &\quad	- f''(\Db+\Dr)(\Db+\Dr)(\ddb+\ddv)\gj\\
		&= - f(\Db+\Dr)(\ddb+\ddv)\gj\\
		 &\quad +f'(\Db+\Dr) T_2\gj,
\end{align*}
where
\begin{equation}\label{eq:T2}
	T_2 = \Eb + [\dda,\Db] + \Er + (k+h)(\ddb+\ddv).
\end{equation}

Hence we complete the proof by showing
\begin{equation}\label{a1}
	T_2 \gj = 0.\tag{{\Romannumeral 2}'}
\end{equation}
Again this can be proved  by computing commutation relations.
We compute
\begin{align*}
	\Eb \Da^{i+1}\jjj &= E^{\ol a} D^{b}D^{c_1}\dots D^{c_i}
		j\langle
		  \nabla^{i+2}\phi,
		  e_{c_i}\dots e_{c_1} e_b e_{\ol a}
		\rangle\\
	&= B_0 D^{c_1}\dots D^{c_i}
		j\langle
		  \nabla^{i+2}\phi,
		  e_{c_i}\dots e_{c_1} e_a e_{\ol a}
		\rangle\\
	&\quad +
		D^{b} E^{\ol a} D^{c_1}\dots D^{c_i}
		j\langle
		  R\nabla^{i}\phi,
		  e_{c_i}\dots e_{c_1} e_b e_{\ol a}
		\rangle
		+ \Da \Eb \Da^{i}\jjj.
\end{align*}
Hence we obtain
\begin{align*}
	[\Eb,\Da] \Da^{i}\jjj &=
	  B_0 D^{c_1}\dots D^{c_i}
		j\langle
		  \nabla^{i+2}\phi,
		  e_{c_i}\dots e_{c_1} e_a e_{\ol a}
		\rangle\\
	&\quad +
		D^{b} E^{\ol a} D^{c_1}\dots D^{c_i}
		j\langle
		  R\nabla^{i}\phi,
		  e_{c_i}\dots e_{c_1} e_b e_{\ol a}
		\rangle\\
	&=
	  2 B_0 D^b D^{c_1}\dots D^{c_{i-1}}
		j\langle
		  R \nabla^{i}\phi,
		  e_{c_{i-1}}\dots e_{c_1} e_a
			 e_b  e_{\ol a}
		\rangle\\
	&\quad +
		D^{b} D^{d} E^{\ol a} D^{c_1}\dots D^{c_{i-1}}
		j\langle
		  - R^e_d \nabla_{e_e} \nabla^{i-1}\phi,
		  e_{c_{i-1}}\dots e_{c_1} e_b e_{\ol a}
		\rangle\\
	&\quad +
		\Da [\Eb,\Da]\Da^{i-1}\jjj.
\end{align*}
Moreover we compute
\begin{align*}
	[[\Eb,\Da],\Da] \Da^{i}\jjj
	 &=
	  2 B_0 D^b D^{c_1}\dots D^{c_{i}}
		j\langle
		  R \nabla^{i+1}\phi,
		  e_{c_{i}}\dots e_{c_1} e_a
			 e_b  e_{\ol a}
		\rangle\\
	&\quad
	  + 2 \Curv bade  D^b D^d E^{\ol a}
	    D^{c_1}\dots D^{c_i}
		j\langle
		  \nabla^{i+1}\phi,
		  e_{c_{i}}\dots e_{c_1}  e_e
		\rangle\\
	&=
	  2 B_0 D^b D^d D^{c_1}\dots D^{c_{i-1}}
		j\langle
		  - R^e_d \nabla_{e_e} \nabla^{i}\phi,
		  e_{c_{i-1}}\dots e_{c_1} e_a
			 e_b  e_{\ol a}
		\rangle\\
	&\quad
	  + 2 \Curv bade  D^b D^d B_0
	    D^{c_1}\dots D^{c_{i-1}}
		j\langle
		  \nabla^{i+1}\phi,
		  e_{c_{i-1}}\dots e_{c_1}
		    e_a  e_e
		\rangle\\
	&\quad
	  + \Da [[\Eb,\Da],\Da]\Da^{i-1}\jjj.
\end{align*}
Hence we obtain
\begin{align*}
	[[[\Eb,\Da],\Da],\Da] \Da^{i}\jjj
	 =
	 6 B_0 \Curv bade D^b D^d
	    D^{c_1}\dots D^{c_i}
		j\langle
		  \nabla^{i+2}\phi,
		  e_{c_{i}}\dots e_{c_1}
			 e_a e_e
		\rangle.
\end{align*}
Thus we have
$$
	[[[[\Eb,\Da],\Da],\Da],\Da] \Da^{i}\jjj = 0.
 $$
Therefore we obtain
$$
	[\Eb,g(\Da)]\jjj = (g'(\Da)S_1 + g''(\Da)S_2
		 + g'''(\Da) B_0 S_3)(\phi),
 $$
where
\begin{align}
	S_1(\phi) &= B_0
		j\langle
		  \nabla^2\phi,
		    e_a  e_{\ol a}
		\rangle
	 + D^a E^{\ol b}
		j\langle
		  R\phi,
		    e_a  e_{\ol b}
		\rangle,
	  \label{eq:S_1}\\
	S_2(\phi) &=
	   B_0 D^a
		j\langle
		  R \nabla\phi,
		  e_b  e_a  e_{\ol b}
		\rangle
	  + \Curv abcd  D^a D^c E^{\ol b}
		  j(\nabla_{e_d}\phi),
	  \label{eq:S_2}\\
	S_3(\phi) &=
	 \Curv abcd  D^a D^c
		j\langle
		  \nabla^{2}\phi,
		  e_b e_d
		\rangle.
	  \label{eq:S_3}
\end{align}

Define
\begin{align}
	C_2 &= \Ric ab D^a \La b,
		\label{eq:C2}\\
	T_3 &= \Curv abcd e^{\ol b} D^a D^c \La d,
		\label{eq:T3}\\
	T_4 &= -2 \Curv abcd e^a D^{\ol b} D^c \La d.
		\label{eq:T4}
\end{align}
Then, in the same way, we compute
\begin{align*}
	[\Er,g(\Da)]\jjj &= g'(\Da)(\Da B_2 - \Dr \dda - B_0 C_2)\jjj,\\
	[\ddb,g(\Da)]\jjj &= (g'(\Da) [\ddb, \Da] + g''(\Da) T_3)\jjj,\\
	[\ddv, g(\Da)]\jjj &= 0,\\
	[[\dda,\Db], g(\Da)]\jjj &= g'(\Da) T_4 \jjj.
\end{align*}

Therefore we obtain
\begin{align*}
	[T_2, g(\Da) ]\jjj
	 &=
		g'(\Da)(S_1  + (k+h)[\ddb,\Da]j
			- B_0 (2C_2 - r\Da)j)(\phi)\\
		&\quad + g''(\Da)(S_2 + (k+h)T_3j)(\phi)
		 + g'''(\Da)B_0S_3(\phi),
\end{align*}
where we use
$$
	T_4 = \Dr\dda - B_0 C_2 - \Da (B_2 - r B_0).
 $$

We have
$$
	S_3 (\phi) = \Da(2 C_2 - r \Da)\jjj.
 $$
Thus, by assuming the relations:
\begin{gather*} \tag{{\Romannumeral 2}''a}\label{y_eight}
	( S_2+(k+h)T_3j -(n-k+h+1)(2C_2 - r \Da)B_0j)(\phi)
	 = - \Da T_2\jjj,\\
 \label{y_seven}\tag{{\Romannumeral 2}''b}
	\left(S_1 + (k+h)[\ddb,\Da]j \right)(\phi)
	 = - (n-k+h)T_2\jjj,
\end{gather*}
we have
\begin{align*}
	[T_2,g(\Da)]\jjj = - g(\Da) T_2\jjj.
\end{align*}
Therefore we have \eqref{a1} and the theorem follows.

To prove \eqref{y_eight} and \eqref{y_seven},
we compute $T_2\jjj$ first. By using Lemma \ref{lem:act_e} and Lemma
\ref{lem:com_ed}, we compute
\begin{align*}
	\Eb\jjj &= - (k+h)\ddb\jjj,\\
	\Er\jjj &= (h (
	  S_4 + S_5) + (k+h)B_1j
	  + (n-k+h)B_2j
	  - (n+1)r B_0j)(\phi),
\end{align*}
where
\begin{align}
	S_4(\phi) &=
	  \Ric ab D^a e^{\ol c}
	    \phi_{b\ol{c}\ol{I}_1,\dots,\ol{I}_h}
	    s^{\ol{I}_1,\dots,\ol{I}_h},
		\label{eq:S_4}\\
	S_5(\phi) &=
	  \Ric ab D^{\ol b} e^c
	    \phi_{\ol a c\ol{I}_1,\dots,\ol{I}_h}
	    s^{\ol{I}_1,\dots,\ol{I}_h}.
		\label{eq:S_5}
 \end{align}
We compute
\begin{align*}
	[\dda,\Db]\jjj &=
	  h(- S_5(\phi) + B_0 j(S_0\phi))  - (n-2k) B_2\jjj\\
	&\quad
	  + h(n-k)r B_0\jjj + B_0 j(\El \phi),
 \end{align*}
where  we use \eqref{eq:spin_curv} and
\begin{align*}
	\Ric ab e^a D^{\ol c} \phi_{\ol cb\ol{I}_1,\dots,\ol{I}_h}
	s^\iii
	 = - B_2\jjj.
 \end{align*}
Hence, by \eqref{eq:ddv2B}, we have
\begin{align*}
	T_2\jjj
	  = h (
		S_4(\phi)
	      + B_0 j(S_0\phi))
	  + (h-1)(n-k+1) r B_0\jjj
	 + B_0 j(\El \phi).
\end{align*}

In order to prove \eqref{y_eight},
we compute first
\begin{align*}
	\Curv abcd   D^a D^c E^{\ol b}
		  j(\nabla_{e_d}\phi)
	 + (k+h)T_3\jjj
	&= - h \Curv abcd D^aD^c e^{\ol e}
	 (\nabla_{e_d}\phi)_{b\ol{e}\ol{I}_1,\dots,,\ol{I}_h}
		s^{\ol{I}_1,\dots,\ol{I}_h}\\
	&= - h( \Da S_4 - B_0C_2j + r \Da B_0j)(\phi),
\end{align*}
where we use
\begin{align*}
	\Ric ad D^a D^b e^{\ol e}
	  (\nabla_{e_d}\phi)_{b\ol e\ol{I}_1,\dots,\ol{I}_h}
		s^{\ol{I}_1,\dots,\ol{I}_h}
	 = - B_0 C_2\jjj.
 \end{align*}
Next we compute using \eqref{eq:spin_curv}
\begin{align*}
	D^a
	  j\langle
	    R \nabla\phi,
	    e_b  e_a  e_{\ol b}
	  \rangle
	&= D^a ( -
	  j\langle
	    R^c_b \nabla_{e_c} \phi,
	     e_a e_{\ol b}
	  \rangle
	+
	  j\langle
	   \nabla_{e_b} R \phi,
	     e_a e_{\ol b}
	  \rangle
	)\\
	&=
	   ( h + 2(n-k+1)) C_2\jjj
		 - h \Da j(S_0\phi)\\
	&\quad
		 - h(n-k+1) r \Da \jjj
		 - \Da j(\El\phi),
\end{align*}
where we use
\begin{align*}
	\Ric ab D^a(\nabla_{e_c}\phi)_{\ol c b\ol{I}_1,\dots,\ol{I}_h}
	s^{\ol{I}_1,\dots,\ol{I}_h}
	 &= 0,\\
	\Ric ab D^c(\nabla_{e_b}\phi)_{\ol ac\ol{I}_1,\dots,\ol{I}_h}
	s^{\ol{I}_1,\dots,\ol{I}_h}
	 &= 0.
 \end{align*}
Therefore we obtain \eqref{y_eight}.

 As for \eqref{y_seven}, we have
\begin{align*}
	D^a E^{\ol b}
		j\langle
		  R\phi,
		    e_a  \ol{e_b}
		\rangle
	+ (k+h) [\ddb, \Da]\jjj = - h S_6(\phi),
\end{align*}
where
\begin{equation}\label{eq:S_6}
	S_6(\phi) = D^a e^{\ol b}
	  \langle
		R\phi, e_a e_{\ol c}
	  \rangle_{c\ol b \ol{I}_1,\dots,\ol{I}_h}
	  s^{\ol{I}_1,\dots,\ol{I}_h}.
\end{equation}
We compute using \eqref{eq:spin_curv}
\begin{align*}
	S_6(\phi)
	 &=
	 (2\Curv abcd D^a e^{\ol e}
	    \phi_{\ol cdb\ol e\ol{I}_1,\dots,\ol{I}_h}
	 +2(h-1) \Curv abcd D^ae^{\ol e}
	    \phi_{b\ol e\ol{I}_1,\ol cd\ol{I}_2,\dots,\ol{I}_h})
	  s^\iii\\
	 & \quad +(h(k+1) + n - 2k) S_4(\phi)
	  + B_0 j(\El\phi).
 \end{align*}
By using \eqref{eq:plucker} and Definition \ref{defn:bk}, we compute
\begin{align*}
	2\Curv abcd D^a e^{\ol e}
	    \phi_{\ol cdb\ol e\ol{I}_1,\dots,\ol{I}_h}
	s^\iii
	= 0
 \end{align*}
and
\begin{align*}
	2 \Curv abcd D^ae^{\ol e}
	    \phi_{b\ol e\ol{I}_1,\ol cd\ol{I}_2,\dots,\ol{I}_h}
	s^\iii
	=
	  -k S_4(\phi)
	 + B_0 j(S_0\phi)
	 + (n-k+1)r B_0\jjj.
 \end{align*}
Thus we have
\begin{align*}
	S_6(\phi)
	= (n-k+h) S_4(\phi)
	 + (h-1)B_0\left(
	  j(S_0\phi)
	  + (n-k+1)r \jjj
	\right)
	+ B_0j(\El\phi).
 \end{align*}
Hence, by using \eqref{eq:lap}, we have \eqref{y_seven}.
\end{proof}

%
\section*{Appendix: List of operators in \S\ref{sec:inv}}\label{sec:apdx}

\newcommand{\apb}{\displaybreak[0]}
\newcommand{\myitem}[1]{\intertext{\indent$\bullet$ #1}}

\begin{align*}
\myitem{Tensorial operators}
	\FF^a_b &= - \FF(\E^a_b) =
		- z^{a\ol \imath I}z^{\ol bjI}\partdel{}{w_{ij}}
	  \reftag{eq:vf}\apb\\
	\FF^{\ol a}_{\ol b} &= \overline{\FF^a_b}\apb\\
	T_0 &= - i(\FF^b_a)R^a_b - i(\FF^{\ol b}_{\ol a})R^{\ol a}_{\ol b}
		\tag{\ref{eq:T0}}\apb\\
	\La a &= L_{e_a} - \omega^b_a i(e_b)
		- [i(e_a),T_0]\reftag{eq:La}\apb\\
	\Lb a &= L_{e_{\ol a}} -
		 \omega^{\ol b}_{\ol a} i(e_{\ol b})
		- [i(e_{\ol a}),T_0]\reftag{eq:Lb}\apb\\
	D^{a} &= e^b i(\FF^{\ol b}_{\ol a})\reftag{eq:D^a}\apb\\
	D^{\ol{a}} &= - e^{\ol b} i(\FF^{\ol a}_{\ol b})\reftag{eq:D^b}\apb\\
	E^{a} &= [\db, D^{a}] + \omega^a_bD^b
	  \reftag{eq:E^b}\apb\\
	E^{\ol a} &= [\db, D^{\ol a}]
	   + \omega^{\ol a}_{\ol b}D^{\ol b}
	  \reftag{eq:E^a}\apb\\
\myitem{Global operators which are used for construction of the inverse
Penrose transform}
	\Da &=  D^{a} \La a
		\reftag{eq:Da}\apb\\
	\Db &= D^{\ol{a}}\Lb a
		\reftag{eq:Db}\apb\\
	\Dr &= - \Ric{a}{b}D^{\ol{b}}D^a
		\reftag{eq:Dr}\apb\\
%
\myitem{Global indecomposable operators}
	B_0 &= e^a D^{\ol a} = D^{a} e^{\ol a}
		\reftag{eq:B0}\apb\\
	B_1 &=  \Ric ab D^a e^{\ol b}
		\reftag{eq:B1}\apb\\
	B_2 &=  \Ric ab e^a D^{\ol b}
		\reftag{eq:B2}\apb\\
	C_1 &= \Ric ab D^{\ol b} \Lb a
		\reftag{eq:C1}\apb\\
	C_2 &= \Ric ab D^a \La b
		\reftag{eq:C2}\apb\\
	\dda &= e^a \La a
		\reftag{eq:da}\apb\\
	\ddb &= e^{\ol a} \Lb a
		\reftag{eq:db}\apb\\
	\Ea &= E^{a}\La a
		\reftag{eq:Ea}\apb\\
	\Eb &= E^{\ol{a}} \Lb a
		\reftag{eq:Eb}\apb\\
\myitem{The other global operators}
	\ddv &= - i(\FF^{\ol b}_{\ol a})R^{\ol a}_{\ol b}
		\tag{\ref{eq:dv}}\apb\\
	\Er &= [\db, \Dr]
		\reftag{eq:Er}\apb\\
	T_1 &= \Curv abcd e^a D^{\ol b} D^{\ol d} \Lb c
		\reftag{eq:T1}\apb\\
	T_2 &= \Eb + [\dda,\Db] + \Er + (k+h)(\ddb+\ddv)
		\reftag{eq:T2}\apb\\
	T_3 &= \Curv abcd e^{\ol b} D^a D^c \La d
		\reftag{eq:T3}\apb\\
	T_4 &= -2 \Curv abcd e^a D^{\ol b} D^c \La d
		\reftag{eq:T4}\apb\\
%
\myitem{Global sections}
	S_0(\phi) &=
	  \Ric ab
	    \phi_{b\ol a\ol{I}_1,\dots,\ol{I}_h} e^{\ol{I}_1,\dots,\ol{I}_h}
	  \reftag{eq:S0}\apb\\
	S_1(\phi) &= B_0
		j\langle
		  \nabla^2\phi,
		    e_a  e_{\ol a}
		\rangle
	 + D^a E^{\ol b}
		j\langle
		  R\phi,
		    e_a  e_{\ol b}
		\rangle
	  \reftag{eq:S_1}\apb\\
	S_2(\phi) &=
	   B_0 D^a
		j\langle
		  R \nabla\phi,
		  e_b  e_a  e_{\ol b}
		\rangle
	  + \Curv abcd   D^a D^c E^{\ol b}
		  j(\nabla_{e_d}\phi)
	  \reftag{eq:S_2}\apb\\
	S_3(\phi) &=
	 \Curv abcd  D^a D^c
		j\langle
		  \nabla^{2}\phi,
		  e_b e_d
		\rangle
	  \reftag{eq:S_3}\apb\\
	S_4(\phi) &=
	  \Ric ab D^a e^{\ol c}
	    \phi_{b\ol{c}\ol{I}_1,\dots,\ol{I}_h} s^\iii
	  \reftag{eq:S_4}\apb\\
	S_5(\phi) &=
	  \Ric ab D^{\ol b} e^c
	    \phi_{\ol a c\ol{I}_1,\dots,\ol{I}_h} s^\iii
	  \reftag{eq:S_5}\apb\\
	S_6(\phi) &= D^a e^{\ol b}
	  \langle
		R\phi, e_a e_{\ol c}
	  \rangle_{c\ol b \ol{I}_1,\dots,\ol{I}_h} s^\iii
	  \reftag{eq:S_6}
\end{align*}

\end{document}